\documentclass[11pt,a4paper]{article}
\usepackage[numbers,sort&compress]{natbib}
\usepackage{moreverb}
\usepackage{appendix}
\usepackage{subcaption} 
\usepackage{graphicx,caption}
\usepackage{epstopdf}
\usepackage{amsmath,amssymb,amsthm}
\usepackage{algorithm}
\usepackage{algpseudocode}
\usepackage{multirow}
\usepackage[colorlinks,bookmarksopen,bookmarksnumbered,citecolor=red,urlcolor=red]{hyperref}
\usepackage{rotating}
\usepackage[text={170mm,240mm},centering]{geometry}
\usepackage[Symbol]{upgreek}
\usepackage{color}
\usepackage{enumerate}
\usepackage{authblk}
\usepackage{booktabs}
\usepackage{threeparttable}
\usepackage{diagbox}
\definecolor{mycolor}{rgb}{0.0,0.0, 1.0}

\newtheoremstyle{temp}
{6pt plus 2pt minus 2pt}
{6pt plus 2pt minus 2pt}
{}
{}
{\itshape}
{}
{\newline}
{}

\theoremstyle{temp}

\begin{document}
	\title{Perturbation Function Iteration Method: A New Framework for Solving Periodic Solutions of Non-linear and Non-smooth Systems}
    \author[2]{LiMin Cao}
	\author[2]{Yanmao Chen}
    \author[1]{Loic Salles}
    \author[2]{Li Wang}
	\author[1]{Zechang Zheng\thanks{Correspondence to: Zechang Zheng, Email: zzheng@uliege.be}$^,$}	
    \affil[1]{Laboratory Vibration of Turbomachines, University of Liege, Belgium}
	\affil[2]{Department of applied mechanics and engineering, Sun Yat-sen University, Shenzhen, P.R. China}
	\date{}
	\maketitle

\begin{abstract}
	Computing accurate periodic responses in strongly nonlinear or even non-smooth vibration systems remains a fundamental challenge in nonlinear dynamics. Existing numerical methods—most notably the Harmonic Balance Method (HBM) and the Shooting Method (SM)—have achieved notable success but face intrinsic limitations when applied to complex, high-dimensional, or non-smooth systems. A key bottleneck is the construction of Jacobian matrices for the associated algebraic equations; while numerical approximations can avoid explicit analytical derivation, they become unreliable and computationally expensive for large-scale or non-smooth problems. To overcome these challenges, this study proposes the Perturbation Function Iteration Method (PFIM), a novel framework built upon perturbation theory. PFIM transforms nonlinear equations into time-varying linear systems and solves their periodic responses via a piecewise constant approximation scheme. Unlike HBM, PFIM avoids the trade-off between Fourier truncation errors and the Gibbs phenomenon in non-smooth problems by employing a basis-free iterative formulation, while significantly simplifying the Jacobian computation. Extensive numerical studies—including self-excited systems, parameter continuation, systems with varying smoothness, and high-dimensional finite element models—demonstrate that PFIM achieves quadratic convergence in smooth systems and maintains robust linear convergence in highly non-smooth cases. Moreover, comparative analyses show that, for high-dimensional non-smooth systems, PFIM attains solutions of comparable accuracy with computational costs up to two orders of magnitude lower than HBM. These results indicate that PFIM provides a robust and efficient alternative for periodic response analysis in complex nonlinear dynamical systems, with strong potential for practical engineering applications.
\end{abstract}

\section{Introduction}
\label{Section1}

    Vibration theory provides a fundamental framework for describing the motion of mechanical and structural systems. In the classical setting, linear vibration theory relies on the principle of superposition, which enables the development of elegant and powerful analytical tools such as modal analysis\cite{TERMEULEN2025111822}, frequency response analysis, and linear control design. These methods have been successfully applied across numerous fields. However, most real-world systems are inherently nonlinear, and their behaviors cannot be captured by linear extrapolations. 
    
    Therefore, nonlinear vibration analysis aims to elucidate how system responses evolve as parameters vary, a task fundamentally addressed through bifurcation analysis. Bifurcations give rise to a wide spectrum of nonlinear phenomena absent in linear theory. For instance, a saddle-node bifurcation \cite{duenas2025saddle} may induce amplitude jumps or solution annihilation, whereas a Hopf bifurcation \cite{mohammed2025hopf} alters equilibrium stability and generates self-sustained limit-cycle oscillations. Period-doubling bifurcations \cite{HU2024115521} can initiate cascades leading to chaos. Even more striking are bifurcations unique to non-smooth systems, such as grazing \cite{wang2025bifurcations} and sliding \cite{10.1115/1.4052882, lu2025experimental} bifurcations, which introduce rich and often abrupt dynamical behaviors.

    The origins of nonlinearities in engineering systems can generally be classified into several categories. Structural nonlinearities arise from elements such as cubic \cite{kuznetsov2025forced} or piecewise-linear \cite{zhang2025dynamic} springs; geometric nonlinearities are associated with large deflections \cite{li2025large}, tension–compression asymmetry \cite{misra2025emergence}, or follower forces \cite{lamy2025sensitivity}; material nonlinearities include plasticity \cite{liu2025tuning}, viscoelasticity \cite{he2025unidirectional}, and nonlinear damping laws \cite{greiner2024model}; and boundary or contact nonlinearities are introduced by clearance \cite{ambrozkiewicz2022influence}, friction \cite{bekesi2025phase}, or impact \cite{10.1016/j.jsv.2025.119150}. For instance, in aeroelastic structures such as aircraft wings, geometric stiffening under large deflection leads to hardening-type nonlinear frequency shifts, while control surface freeplay introduces piecewise stiffness \cite{ZHENG2022103440} that can trigger limit-cycle oscillations and bifurcations.

    Just as the principle of superposition forms the foundation of linear vibration theory, perturbation theory constitutes the core of classical nonlinear vibration analysis. Early computational approaches were predominantly analytical and built upon perturbation expansions. Representative examples include the multiple scale method \cite{salih2014method}, the Lindstedt-Poincaré method \cite{amore2005improved}, and the averaging method \cite{sanders2007averaging}, all of which systematically expand the solution in terms of a small parameter characterizing the system’s nonlinearity. The key idea is to express the response as a power series in this parameter, substitute it into the governing equations, and solve the resulting hierarchy of linear differential equations at successive orders.

    Although these methods have provided valuable insights and analytical solutions for a variety of nonlinear systems, their limitations are well recognized: (i) they are restricted to weakly nonlinear systems, where a small parameter can be clearly identified, and (ii) they are difficult to extend to multi-degree-of-freedom systems or systems with strong nonlinearities, where convergence of the perturbation series may fail.

    In recent decades, research focus has gradually shifted from analyzing equilibria to investigating the periodic responses of nonlinear systems. With the advent of modern high-performance computing, numerical methods have become the dominant approach for studying periodic solutions. Among these, the harmonic balance method (HBM) \cite{10.1016/j.jsv.2024.118925, YAN20231419} and the shooting method (SM) \cite{10.1115/1.4038327, LOGHMAN2022116521} are the most widely used and represent the primary techniques compared in this work.

    To illustrate their principles and highlight their similarities and differences, we begin with a generalized nonlinear dynamical system formulated in state-space form:
    \begin{equation*}
	   \dot{x} = f(x, t),
    \end{equation*}
    where $x$ represents the state vector and the overdot denotes differentiation with respect to time $t$. To compute periodic responses of this system, one must additionally impose the periodic boundary condition.
    \begin{equation*}
	   x(t+T) = x(t),
    \end{equation*}
    where $T$ is the (maybe unknown) period of the solution. On this basis, the fundamental ideas of HBM and SM, as well as their similarities and differences, are carefully examined.

\subsection*{Harmonic Balance Method}

    HBM is a representative basis-function approach for solving periodic responses of nonlinear systems. Its fundamental idea is to approximate the periodic solution $x(t)$ using a truncated Fourier series of $H$ harmonics:
    \begin{equation*}
    x(t) \approx x_H(t) = a_{0} + \sum_{i=1}^{H} [a_{i}\cos(i\omega t) + b_{i}\sin(i\omega t)],
     \end{equation*}
    where $a_{i}$ and $b_{i}$ are the unknown Fourier coefficients collected in the vector $C_{H} = [a_0,a_1,...a_n,b_1,...,b_n]$.  
    Due to the Fourier series representation, the periodic boundary condition $x(t+T)=x(t)$ is satisfied automatically.  

    Substituting $x_H(t)$ into the governing equations yields the residual function $R_H(t)$
      \begin{equation*}
	         R_H(t) =\dot{x}_H - f(x_H,t). 
     \end{equation*}
    To convert the original differential equation into a set of algebraic equations, the Galerkin procedure is applied, exploiting the orthogonality of the Fourier basis functions:
      \begin{equation*}
      G_{H} = \int_{0}^{T} \mathbf{W}_H^{\top}(t)\,R_H(t)\, dt = 0,
      \end{equation*}
    where \(\mathbf{W}_H(t) = [1,\cos(\omega t),\cos(2\omega t),\dots,\cos(H\omega t),\sin(\omega t),\sin(2\omega t),\dots,\sin(H\omega t)]\) is the weight vector.  

\begin{figure*}[!htbp]
    	\centering
    	\raisebox{1pt}{\begin{subfigure}[t]{0.7\textwidth}
    			\centering
    			\includegraphics[width=\linewidth]{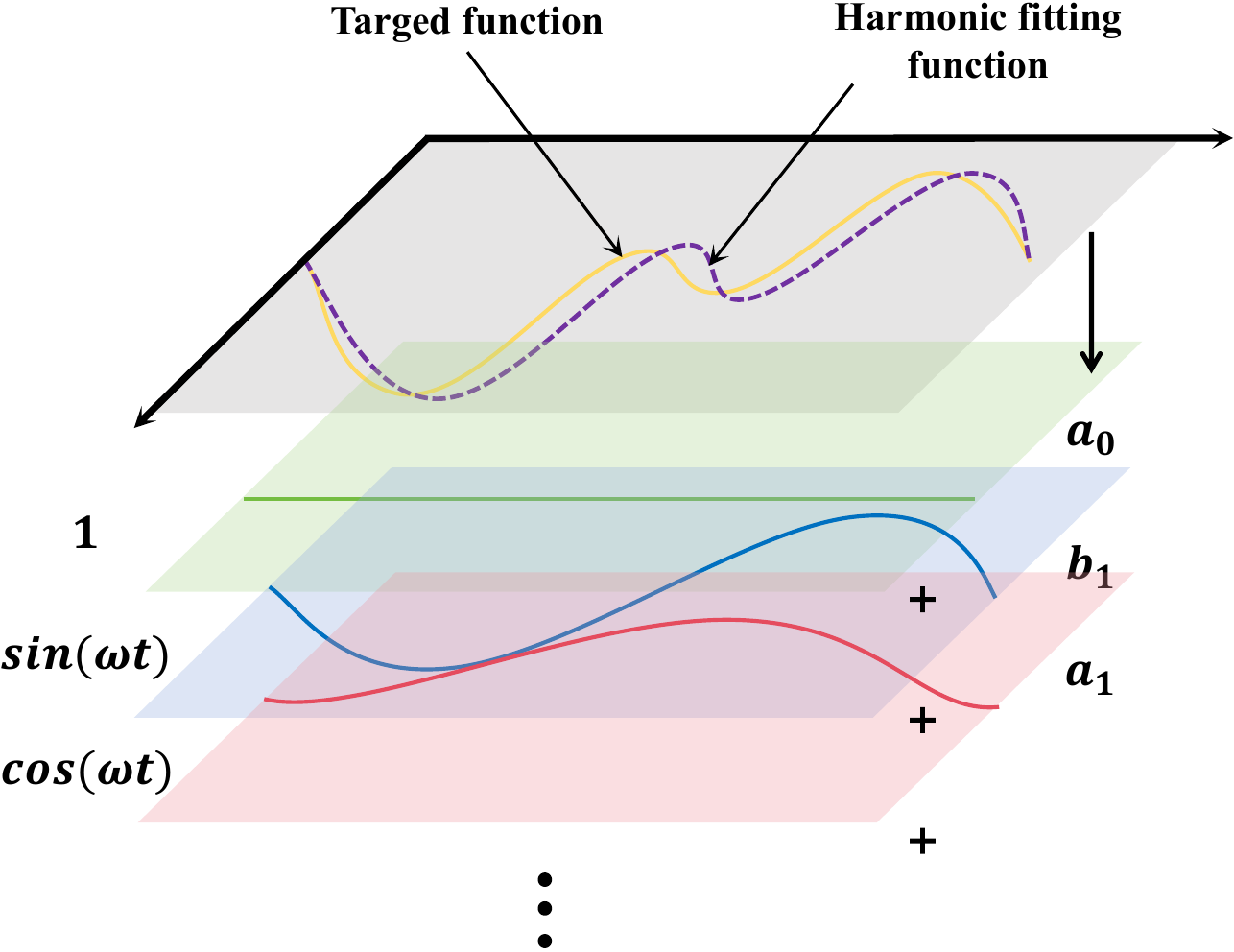}
    	\end{subfigure}}	
    	\caption{Schematic diagram of harmonic balance method.}
    	\label{harmonic_scheme}
    \end{figure*}

    The resulting nonlinear algebraic system $G_{H}=0$ is then solved iteratively, typically using Newton–Raphson iterations, until the residual norm $\|G_H\|$ falls below a prescribed tolerance. A schematic illustration of the solution procedure of HBM is provided in Fig.~\ref{harmonic_scheme} to provide an intuitive overview.

    HBM and its numerous variants have become one of the most classical and widely adopted approaches to analyze the periodic responses of modern nonlinear systems \cite{hall2013harmonic, 10.1115/1.4066216, PEI2022106220}. For smooth systems, representing periodic motions using Fourier series is highly efficient and accurate. However, when applied to non-smooth systems, this very feature turns into a drawback.

    Although several strategies have been proposed, such as the alternate frequency-time (AFT) method \cite{CHEN2023109805} or collocation-based techniques \cite{krack2019harmonic, dai2012simple} - to mitigate some of the computational difficulties, two fundamental limitations remain unresolved. First, the convergence of the Fourier series deteriorates significantly for non-smooth functions \cite{zhou2025enriching}, necessitating much higher harmonic orders than in smooth problems. Second, approximating discontinuous functions with Fourier series inevitably introduces Gibbs oscillations \cite{gottlieb1997gibbs}, and this phenomenon does not vanish even as the number of harmonics tends to infinity. Taken together, these two issues form a fundamental contradiction: the first compels the use of a large number of harmonics to achieve accuracy, while the second guarantees that no matter how many harmonics are included, the approximation error due to Gibbs oscillations cannot be completely eliminated.
    
\subsection*{Shooting Method}

    In contrast to HBM, SM approaches the problem from the time domain. Its core idea is to determine an initial state $x_{0}$ such that the periodic boundary condition is satisfied after numerical integration over one period. The trajectory is obtained using a suitable time-integration scheme—such as the classical fourth-order Runge–Kutta method, the Newmark–$\beta$ method, or other problem-specific algorithms.  

\begin{figure*}[!htbp]
    	\centering
    	\raisebox{1pt}{\begin{subfigure}[t]{0.7\textwidth}
    			\centering
    			\includegraphics[width=\linewidth]{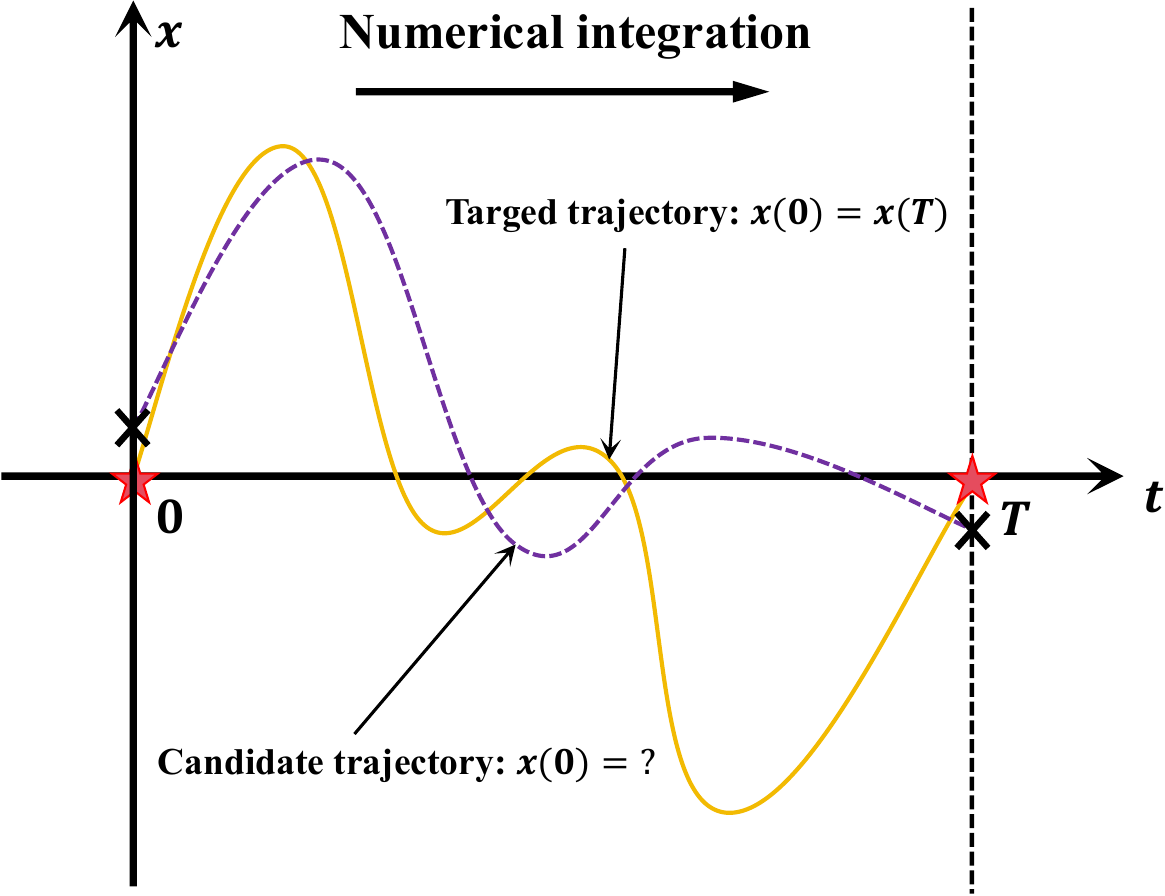}
    	\end{subfigure}}	
    	\caption{Schematic diagram of shooting balance method.}
    	\label{shooting_scheme}
    \end{figure*}

    Let $\psi(x_{0},T)$ denote the state of the system after integrating over one period $T$ starting from $x(0)=x_{0}$:
    \begin{equation*}
    \psi(x_{0},t) = x_{0} + \int_{0}^{t} f(x,s)\, ds,
    \end{equation*}
    where the integral represents the numerical flow map of the governing equations. Substituting this expression into the periodic boundary condition $x(T)=x(0)$ yields the nonlinear algebraic equation
    \begin{equation*}
    G_{S} = \psi(x_{0},T) - x_{0} = 0.
    \end{equation*}
    This system is then solved iteratively—typically using Newton–Raphson iterations—until the norm of $G_{S}$ falls below a prescribed tolerance. For clarity, Fig.~\ref{shooting_scheme} provides a schematic illustration of the SM solution procedure.

    Unlike HBM, SM can be viewed as a time-domain strategy that leverages numerical integration to obtain periodic solutions. One of its key advantages is its ability to accurately capture unstable periodic orbits, since the integration is performed only over a single period and avoids the long-term divergence issues encountered in direct time integration.

    As its core relies on time integration, most improvements to SM have focused on developing more efficient and more accurate integration schemes, such as adaptive Runge–Kutta methods with embedded error control, implicit single-step methods for stiff systems, and symplectic schemes for Hamiltonian dynamics. These advances have significantly enhanced the robustness and computational efficiency of SM.

    Nevertheless, two inherent limitations remain unchanged. One is that, SM requires repeated numerical integration within each Newton iteration, leading to high computational cost \cite{charroyer2018self, lim2024proper}, especially for systems with many degrees of freedom or very long periods. Next, and most critically, the Jacobian of the shooting function $G_S$ is often difficult to obtain. For smooth systems, analytical Jacobians may still be derived but can be cumbersome; for non-smooth systems, deriving an analytical Jacobian is often intractable due to discontinuities and non-differentiable terms. As a result, finite-difference approximations or variational equations are commonly used, but these approaches incur substantial additional computational cost and can suffer from severe numerical instability, especially when discontinuities induce large sensitivity in the system response.

    These factors collectively limit the scalability of SM for large-scale engineering problems, despite its accuracy and versatility.

\subsection*{Motivation and Concept of the Proposed Method}

    Although HBM and SM both reformulate the problem into a set of algebraic equations and solve them iteratively, they follow fundamentally different solution philosophies. HBM begins with an assumed periodic representation of the solution—typically a truncated Fourier series—and adjusts its coefficients until the governing equations are satisfied in a weighted-residual sense. In contrast, SM starts from a trajectory that strictly satisfies the governing differential equations (via time integration) and iteratively adjusts the initial condition $x_0$ until the periodic boundary condition is met.

    These complementary perspectives highlight that both methods ultimately rely on Newton-type iterations in a reduced algebraic space, and both require the Jacobian of the target equations with respect to the chosen variables. However, for non-smooth systems, obtaining accurate Jacobians becomes particularly challenging: analytical Jacobians are often unavailable, while numerical Jacobians incur high computational cost and can be unstable near discontinuities. Moreover, SM is computationally expensive because it performs repeated time integrations, whereas HBM suffers from slow convergence and Gibbs oscillations when approximating non-smooth responses.

    Motivated by these challenges, this study introduces the Perturbation Function Iteration Method (PFIM). Unlike HBM and SM, which rely on auxiliary algebraic variables and finite Fourier representations, PFIM performs Newton updates directly on the governing equations using the Jacobian with respect to the state variables. The solution is treated as a basis-free function rather than a truncated series, and its trajectory is computed through high-accuracy time-integration schemes. This formulation preserves robustness in the presence of non-smooth nonlinearities while ensuring accurate and efficient solution convergence.

    The remainder of this paper is organized as follows: Section 2 presents the theoretical formulation, implementation procedure, and complexity analysis of PFIM. Section 3 demonstrates its application to smooth systems, using the Van der Pol and Duffing oscillators to illustrate its capability in handling unknown frequencies and continuation problems. Section 4 investigates PFIM’s performance in non-smooth systems with varying degrees of continuity ($C^{1}$, $C^{0}$, $C^{-1}$), highlighting its convergence characteristics. Section 5 applies PFIM to a high-dimensional finite-element cantilever beam model, comparing its computational accuracy and efficiency against HBM. Section 6 summarizes the main findings and discusses the implications of PFIM for nonlinear vibration analysis.

\section{Methodology}
\label{Section2}
 \subsection{Basic Idea}
Consider a nonlinear system in state-space form:
\begin{equation}
	\dot{x} = f(x,t),
	\label{eq1}
\end{equation}
where $\dot{x}$ denotes the derivative of $x$ with respect to time $t$. To compute its periodic solution, first extract the oscillation frequency $\omega$ and introduce the dimensionless time $\tau$ scaled over $[0, 2\pi]$. This transforms Eq.~\eqref{eq1} into:
\begin{equation}
	\omega x' = f(x,\tau),
	\label{eq2}
\end{equation}
where $'$ indicates differentiation with respect to $\tau$. The periodic solution $x_*(\tau)$ and true frequency $\omega_*$ satisfies $R(x_*, \omega_*, \tau) \equiv 0$ for $\tau \in [0, 2\pi]$, where $R(x, \omega, \tau) := f(x,\tau) - \omega x'$.

Given an initial periodic guess $x_0(\tau)$ and frequency guess $\omega_0$ such that $R(x_0, \omega_0, \tau) \neq 0$, express the exact solution as $x_*(\tau) = x_0(\tau) + \Delta x(\tau)$ and true frequency $\omega_* = \omega_0 + \Delta \omega$. Substituting this into $R(x_*, \omega_*, \tau) = 0$ and expanding via perturbation yields:
\begin{equation}
	R(x_*, \omega_*, \tau) = f(x_*,\tau) - \omega_* x'_* 
	= \underbrace{f(x_0,\tau) - \omega_0 x'_0}_{R(x_0, \omega_0, \tau)} 
	+ \frac{\partial f}{\partial x}\bigg|_{x_0} \Delta x 
	- \omega_0 \Delta x' - \Delta \omega x'_0 
	+ \mathcal{O} 
	= 0.
	\label{eq3}
\end{equation}
Neglecting higher-order terms results in a linear differential equation for the correction $\Delta x(\tau)$ and $\Delta \omega$:
\begin{equation}
	\Delta x' = \frac{1}{\omega_0}(\frac{\partial f}{\partial x}\bigg|_{x_0} \Delta x + R(x_0,\omega_0,\tau)-\Delta \omega x'_0).
	\label{eq4}
\end{equation}
For simplicity, Eq.~\eqref{eq4} can be rewritten as:
\begin{equation}
	\mu' = Q(\tau) \mu + P(\tau) + F(\tau) \nu,
	\label{eq5}
\end{equation}
where $\mu = \Delta x$, $Q(\tau) = \frac{1}{\omega_0} \left. \frac{\partial f}{\partial x} \right|_{x_0(\tau)}$, $P(\tau) = \frac{1}{\omega_0} R(x_0, \tau)$, $F(\tau)=-\frac{x'_0}{\omega_0}$ and $\nu = \Delta \omega$. According to the theory of linear differential equations, the general solution of Eq.~\eqref{eq5} is
\begin{equation}
	\mu(\tau) = \Phi(\tau, 0) \mu(0) + \Gamma(\tau) + \Pi(\tau) \nu.
	\label{eq6}
\end{equation}
Here, $\Phi(\tau, 0)$ is the state transition matrix for the homogeneous system, given by 
\begin{equation}
	\Phi(\tau, 0) = \exp\left(\int_0^{\tau} Q(s)  ds\right),
\end{equation}
which satisfies the matrix differential equation:
\begin{equation}
	\frac{\partial}{\partial \tau} \Phi(\tau, \tau_0) = Q(\tau) \Phi(\tau, \tau_0).
	\label{eq7}
\end{equation}
The state transition matrix propagates solutions such that $\mu(\tau_1) = \Phi(\tau_1, \tau_2) \mu(\tau_2)$ for any $\tau_1, \tau_2 \in [0, 2\pi]$, with $\Phi(0, 0) = I$ (identity matrix). The particular solutions are defined as:
\begin{equation}
    \left\{
    \begin{aligned}
        \Gamma(\tau) &= \textstyle\int_0^{\tau} \Phi(\tau, s) P(s)  ds \\
        \Pi(\tau) &= \textstyle\int_0^{\tau} \Phi(\tau, s) F(s)  ds
    \end{aligned},
    \right.
\end{equation}
and satisfy the respective differential equations:
\begin{equation}
    \left\{
    \begin{aligned}
        \Gamma'(\tau) &= Q(\tau) \Gamma(\tau) + P(\tau) \\
        \Pi'(\tau) &= Q(\tau) \Pi(\tau) + F(\tau)
    \end{aligned}.
    \right.
    \label{eq8}
\end{equation}

For periodicity, $\mu(\tau)$ must satisfy the boundary condition $\mu(2\pi) = \mu(0)$. Substituting $\tau = 2\pi$ into Eq.~\eqref{eq6} yields:
\begin{equation}
	\mu(2\pi) = \Phi(2\pi, 0) \mu(0) + \Gamma(2\pi) + \Pi(2\pi) \nu.
	\label{eq9}
\end{equation}
In Eq.~\eqref{eq9}, the introduction of $\nu$ results in one more unknown than the number of equations. To resolve this, an additional phase condition must be introduced. Denote this additional phase condition as the following linear equation:

\begin{equation}
	\Upsilon
	\begin{bmatrix}
		\mu(0) \\
		\nu
	\end{bmatrix} 
	= \Xi,
	\label{eq9.5}
\end{equation}
where $\Upsilon$ is an $(n+1)$-dimensional row vector, $\top$ denotes the transpose, and $\Xi$ is a scalar.

For systems with known frequency, such as externally excited systems where the excitation frequency coincides with the target periodic response frequency, we have $\omega_* = \omega_0$. This leads to the additional phase condition $\Delta \omega = 0$. Therefore, the specific forms in Eq.~\eqref{eq9.5} become:
\begin{equation}
	\Upsilon = [O,\ 1], \quad \Xi = 0,
\end{equation}
where $O$ is an $n$-dimensional row vector with all zero elements.

For self-excited systems, $\Delta \omega \neq 0$. Since the additional phase constraint is imposed on the initial point $x(0)$, one may adopt the orthogonal phase condition commonly introduced in the shooting method for autonomous systems \cite{nayfeh2008applied}:
\begin{equation}
	\left\langle \left.\frac{d x_0}{d \tau}\right|_{\tau = 0}, \Delta x(0) \right\rangle = 0.
\end{equation}
This condition is motivated geometrically: if the initial point does not lie on the desired periodic orbit, the correction increment must contain a component orthogonal to the derivative at that point. This phase condition is general and ensures convergence as long as $x(0)$ is sufficiently close to the periodic solution. Thus, the quantities in Eq.~\eqref{eq9.5} take the form:
\begin{equation}
	\Upsilon = \left[\left.\frac{d x_0}{d \tau}\right|_{\tau = 0} ^\top,\ 0\right], \quad \Xi = 0.
\end{equation}

Moreover, the form presented in Eq.~\eqref{eq9.5} can be directly applied in continuation schemes to trace amplitude-frequency curves. When using pseudo-arclength continuation to vary the excitation frequency $\omega$ in externally excited systems, $\omega$ can no longer be treated as a fixed parameter. Instead, both $x(0)$ and $\omega$ are regarded as functions of an arclength parameter $s$, denoted as $x_0(s)$ and $\omega(s)$. After obtaining the solution $(x^i(0), \omega^i)$ at the $i$-th continuation step, the derivatives $\frac{d x^i(0)}{d s}$ and $\frac{d \omega^i}{d s}$ are computed to predict the next solution point:
\begin{equation}
	x^{i+1}_{\text{pre}}(0) = x^i(0) + \frac{d x^i(0)}{d s} \Delta s, \quad \omega^{i+1}_{\text{pre}} = \omega^i + \frac{d \omega^i}{d s} \Delta s,
\end{equation}
where $\Delta s$ is the step size. The prediction is then corrected under the constraint:
\begin{equation}
	N(x^{i+1}_{\text{cor}}(0), \omega^{i+1}_{\text{cor}}) = \left\langle x^{i+1}_{\text{cor}}(0) - x^i(0), \frac{d x^i(0)}{d s} \right\rangle + \left\langle \omega^{i+1}_{\text{cor}} - \omega^i, \frac{d \omega^i}{d s} \right\rangle - \Delta s = 0.
\end{equation}
Correspondingly, the parameters in Eq.~\eqref{eq9.5} are given by:
\begin{equation}
	\Upsilon = \left[\frac{d x^i(0)}{d s} ^\top, \ \frac{d \omega^i}{d s}\right], \quad \Xi = N(x^{i+1}_{\text{cor}}(0), \omega^{i+1}_{\text{cor}}).
\end{equation}

Finally, combining Eq.~\eqref{eq9} with the additional phase condition yields the following algebraic system:
\begin{equation}
	\begin{bmatrix}
		I - \Phi(2\pi) & -\Pi(2\pi) \\
		\multicolumn{2}{c}{\Upsilon}
	\end{bmatrix}
	\begin{bmatrix}
		\mu(0) \\
		\nu
	\end{bmatrix}
	=
	\begin{bmatrix}
		\Gamma(2\pi) \\
		\Xi
	\end{bmatrix}.
	\label{eq10}
\end{equation}
By incorporating the appropriate phase condition according to the system type, one can solve Eq.~\eqref{eq10} for $\Delta x(0)$ and $\Delta \omega$. The periodic correction term $\Delta x(\tau) = \mu(\tau)$ is then obtained via Eq.~\eqref{eq6}.

Due to neglected higher-order terms, $x_0(\tau) + \Delta x(\tau)$ and $\omega_0 +\Delta \omega$ are not yet the exact solution $x_*(\tau)$ and $\omega_*$. However, by updating the initial guess $x_0 \leftarrow x_0 + \Delta x, \omega_0 \leftarrow \omega_0 + \Delta \omega$ and iterating this correction process, the solution converges to $x_*$ and $\omega_*$. This iterative procedure constitutes the PFIM. Through successive corrections $\Delta x(\tau)$ and $\Delta \omega$, the initial approximation progressively converges to the true periodic solution of Eq.~\eqref{eq2}.

\subsection{Piecewise Constant Approximation}
Although the above process can obtain the final periodic solution, the theoretical expressions of the correction terms in each iteration are strictly complex. For low-dimensional systems, it is possible to derive specific expressions of $\Delta x (\tau)$ per iteration. However, for high-dimensional systems, theoretical derivation of each PFIM iteration becomes impractical.

To circumvent solving analytical expressions for high-dimensional linear differential equations, we employ the piecewise constant approximation (PCA) to compute numerical solutions instead of analytical solutions. For linear equations such as Eq.~\eqref{eq6}, PCA method enables high-precision integration while maintaining high computational efficiency. The procedure for combining PFIM with the PCA is presented below.

In the PCA, the entire period is uniformly partitioned into $n$ intervals: $[\tau_0, \tau_1, \dots, \tau_n]$ with $\tau_0=0$ and $\tau_n=2\pi$. The interval length is denoted by $\Delta \tau = \tau_{k+1} - \tau_k$. Within each interval $[\tau_i, \tau_{i+1}]$, the linear time-varying differential equation is approximated as a linear time-invariant system with constant coefficients:
\begin{equation}
	\mu' = Q_i \mu + P_i + F_i \nu,
	\label{eq11}
\end{equation}
where $Q_i = \frac{Q(\tau_i) + Q(\tau_{i+1})}{2}$, $P_i = \frac{P(\tau_i) + P(\tau_{i+1})}{2}$ and $F_i = \frac{F(\tau_i) + F(\tau_{i+1})}{2}$. The general solution in this interval is given by:

\begin{equation}
	\mu_{i+1} = \Phi_i \mu_i + \Gamma_i + \Pi_i \nu,
	\label{eq13}
\end{equation}
in which 
\begin{equation}
	\Phi_i = \exp \left( Q_i \Delta \tau \right),
	\label{eq12_1}
\end{equation}

\begin{equation}
    \left\{
    \begin{aligned}
        \Gamma_i &= \int_{0}^{\Delta \tau} \exp \left( Q_i(\Delta \tau - s) \right)  P_i  ds = \left( \exp \left( Q_i \Delta \tau \right) - I \right) Q_i^{-1} P_i \\
        \Pi_i &= \int_{0}^{\Delta \tau} \exp \left( Q_i(\Delta \tau - s) \right) F_i  ds = \left( \exp \left( Q_i \Delta \tau \right) - I \right) Q_i^{-1} F_i
    \end{aligned}.
    \right.
    \label{eq12_2_3}
\end{equation}

Through recursive relations, the expressions for the fundamental and particular solutions after one full period can be obtained as:

\begin{equation}
    \left\{
    \begin{aligned}
        \Phi_{\text{total}} &= \Phi_{n-1} \Phi_{n-2} \cdots \Phi_{0} \\
        \Gamma_{\text{total}} &= \Gamma_{n-1} + \Phi_{n-1} \Gamma_{n-2} + \cdots + \Phi_{n-1} \cdots \Phi_{1}  \Gamma_{0} \\
        \Pi_{\text{total}} &= \Pi_{n-1} + \Phi_{n-1} \Pi_{n-2} + \cdots + \Phi_{n-1} \cdots \Phi_{1}  \Pi_{0}
    \end{aligned}.
    \right.
    \label{eq14}
\end{equation}

Substituting these back into Eq.~\eqref{eq9} yields:
\begin{equation}
	\mu_n = \Phi_{\text{total}} \mu_0 + \Gamma_{\text{total}} + \Pi_{\text{total}} \nu.
	\label{eq15}
\end{equation}

This approach eliminates the need to compute analytical expressions for $\Delta x(\tau)$ after each iteration. Instead, it only requires computing $\Delta x(\tau_i)$ at discrete time points $i=0,1,2,\ldots, n$ and updating $x_{0}(\tau_i)$ accordingly. The advantage of PCA over conventional numerical integration methods lies in its fundamental approach: traditional methods are inherently approximate. They rely on calculating and combining slopes at multiple points to predict subsequent values, a process that introduces inevitable truncation errors. Although higher-order methods can reduce the error per step, the error persists and accumulates over numerous steps. 

In contrast, PCA leverages an exact theoretical result for linear systems—the matrix exponential. For constant linear differential equations, the matrix exponential provides an exact solution, completely free from truncation error. Consequently, PCA achieves higher accuracy than traditional numerical integration schemes. Furthermore, the number of integration steps required by traditional methods can become exceedingly large for certain types of systems, such as high-dimensional or stiff systems. PCA, however, is not similarly affected by these system properties and typically requires far fewer total steps. Therefore, employing PCA as the numerical solver for PFIM is undoubtedly the most reliable approach. The pseudo code of PFIM with PCA will be demonstrated following
\begin{algorithm}
	\caption{The pseudo code of PFIM with PCA}
	\begin{algorithmic}[1]
		\State \textbf{Step 1:} Determine the initial state $x_0(\tau_i), i=0,1,\dots n$ and initial frequency $\omega_0$. 
		\State \hspace{\algorithmicindent} Give absolute tolerance $tol_{a}$ and relative tolerance $tol_{r}$.
		\State \textbf{Step 2:} Calculate $R(x_0(\tau_i), \omega_0,\tau_i)$ on each $\tau_i$. 
		\State \hspace{\algorithmicindent} Compute average error $e_{a} = \|\frac{\sum_{i=0}^{n}R(x_0(\tau_i), \omega_0, \tau_i)}{n+1}\|$. 
		\State \hspace{\algorithmicindent} If iteration $>1$, calculate $e_{r} = \frac{1}{n+1} \sum_{i=0}^{n} \frac{\|\Delta x(\tau_i)\|}{|x_0(\tau_i)|}$. \label{step2}
		\While{true} \label{loopStart} 
		\State \textbf{Step 3:} Check if $e_a < tol_a$ or $e_r < tol_r$
		\If{$e_a < tol_a$ \textbf{or} $e_r < tol_r$} 
		\State \textbf{break} \Comment{Exit loop if either condition satisfied}
		\Else
		\State \textbf{Step 4:} Compute $Q_i$, $P_i$, $F_i$ $\Phi_i$, $\Gamma_i$ and $\Pi_i$ for each interval.
		\State \hspace{\algorithmicindent} Calculate $\Phi_{\text{total}}$, $\Gamma_{\text{total}}$, $\Pi_{\text{total}}$ via recurrence.
		\State \textbf{Step 5:} Solve $\Delta x(\tau_0)$ using periodic BCs with additional phase condition.
		\State \hspace{\algorithmicindent} Compute remaining $\Delta x(\tau_i)$ and $\Delta \omega$ with $\Phi_i$, $\Gamma_i$ and $\Pi_i$.
		\State \textbf{Step 6:} Update: $x_0(\tau_i) \gets x_0(\tau_i) + \Delta x(\tau_i)$, $\omega_0 \gets \omega_0 + \Delta \omega$ \Comment{Return to Step 2}
		\EndIf
		\EndWhile
		\State Output final solution $x_0(\tau_i)$ and $\omega_0$
	\end{algorithmic}
\end{algorithm}

Compared with HBM, both PFIM and HBM start from an initial guess solution and use Newton iteration to gradually approach the exact solution. However, PFIM does not restrict the solution function to any specific basis representation. Instead, it employs a numerical function as the iterative entity, eliminating the constraints imposed by basis functions on the solution representation. This advantage becomes particularly pronounced when dealing with non-smooth problems.

In contrast to SM, which first computes the solution of a nonlinear system over one period starting from an initial point (satisfying the original equations but not the periodic boundary conditions) and then corrects the initial value through Newton iteration on the boundary-condition-derived objective equation, PFIM adopts a different approach: it postulates the existence of a periodic solution that does not satisfy the original equations, then iteratively refines it until sufficient convergence to the exact solution is achieved. Tab.~\ref{tab_3_3_method} summarizes the distinctive characteristics of PFIM, HBM, and SM.
\begin{table}[!ht]
	\centering
	\caption{The main characteristics of PFIM, HBM and SM solutions}
	\label{tab_3_3_method}
	\begin{tabular}{cccc}
		\toprule
		Method & Type of Solution & Periodic Boundary Condition& Equation Residual \\\midrule
		PFIM& 	Numerical solution  &Strongly satisfied& Satisfied in a weak sense \\
		HBM& Fourier series solution & Strongly satisfied & Satisfied in a weak sense\\
		SM& Numerical solution & Satisfied in a weak sense & Strongly satisfied\\
		\bottomrule
	\end{tabular}
\end{table}

\subsection{Computational Complexity Analysis}
Furthermore, the computational complexity per iteration of the three methods—PFIM, HBM, and SM—is analyzed below to illustrate the computational efficiency of PFIM.

The computational complexity of a single iteration of PFIM is first discussed. Assume the system dimension is \(N\) and the number of integration steps in PCA is \(n_p\). The primary computational cost in PFIM lies in computing the matrix exponential, i.e., \(\Phi_i = \exp(Q_i \Delta \tau)\). Typically, the matrix exponential is approximated using Taylor series expansion:
\begin{equation*}
	\exp(A) \approx \sum_{m=0}^{M} \frac{A^m}{m!},
\end{equation*}
where \(M\) is the truncation order of the Taylor expansion. For an \(N \times N\) matrix, the complexity of one matrix multiplication is \(\mathcal{O}(N^3)\). Hence, the complexity of computing one matrix exponential is \(\mathcal{O}(M N^3)\). With \(n_p\) integration steps, the total complexity becomes \(\mathcal{O}(n_p M N^3)\). Moreover, since the step size in the PCA is generally small, the truncation order \(M\) can also be kept small. Therefore, the overall computational cost of PFIM primarily depends on \(n_p\) and the system dimension \(N\), yielding an overall complexity of \(\mathcal{O}(n_p N^3)\).

When considering the computational complexity of HBM, let $H$ represent the number of harmonic terms retained. The main computational costs arise from  solving the linear system \(J \Delta C_H = R\). The dimension of \(J\) is \(NH \times NH\) since each degree of freedom has \(H\) harmonic coefficients. Solving the linear system incurs a computational cost of \(\mathcal{O}((NH)^3) = \mathcal{O}(H^3 N^3)\).

For SM, the primary computational cost lies in computing the Jacobian matrix \(\Psi\). Denote the number of integration steps as \(n_s\). The governing equation for \(\Psi\) is \(\dot{\Psi} = J_f(t, x) \Psi\), where \(J_f(t, x) = \partial f / \partial x\). Since this is a linear system, the cost per integration step is \(\mathcal{O}(N^2)\). As \(\Psi\) has \(N\) columns, the cost per step for the full Jacobian is \(\mathcal{O}(N^3)\). With \(n_s\) integration steps, the overall complexity becomes \(\mathcal{O}(n_s N^3)\).

Tab.~\ref{tab_2} summarizes the major computational steps and corresponding complexities for the three methods. It can be observed that all methods scale cubically with the system dimension \(N\), implying that the influence of \(N\) is similar across methods. The key factors affecting computational efficiency are the method-specific parameters: \(n_p\) for PFIM, \(H\) for HBM, and \(n_s\) for SM.

\begin{table}[H]
	\centering
	\caption{Computational complexity comparison of periodic solution methods}
	\label{tab_2}
	\begin{tabular}{lcc}
		\toprule
		Method & Dominant Step & Computational Complexity \\
		\midrule
		PFIM & Matrix exponential & $\mathcal{O}(n_p N^3)$ \\
		HBM & Linear system solution & $\mathcal{O}(H^3 N^3)$ \\
		SM & Jacobian computation & $\mathcal{O}(n_s N^3)$ \\
		\bottomrule
	\end{tabular}
\end{table}

In smooth systems, HBM often achieves satisfactory efficiency since a small \(H\) is sufficient for convergence. However, in non-smooth systems, a large \(H\) is generally required, leading to a sharp increase in computational cost due to the cubic dependence on \(H\). Although SM scales linearly with \(n_s\), the required number of steps \(n_s\) can become very large in high-dimensional systems to maintain accuracy and numerical stability. While no universally accepted relationship between \(n_s\) and system dimension exists, numerous studies \cite{charroyer2018self, lim2024proper} on shooting methods indicate that the cost of Jacobian computation increases dramatically with system size, indirectly suggesting that \(n_s\) must also increase significantly.

In contrast, PFIM uses the PCA to integrate the linear differential equations, requiring far fewer steps (\(n_p\)) than SM (\(n_s\)). Moreover, \(n_p\) is largely independent of system dimension, and non-smoothness does not significantly increase \(n_p\). These properties underline the computational advantage of PFIM over both HBM and SM.

It is also worth noting that PFIM holds potential for further acceleration. On one hand, the main cost in PFIM lies in computing matrix exponentials. For very high-dimensional systems (e.g., \(N > 100\)), model order reduction techniques—such as Krylov subspace approximation or SVD-based reduction—can reduce the cost of a single matrix exponential from \(\mathcal{O}(N^3)\) to \(\mathcal{O}(n^3)\), where \(n \ll N\). This offers PFIM an additional efficiency advantage over SM in large-scale systems. On the other hand, each integration step in PFIM is independent, enabling parallel computation. In contrast, SM requires sequential integration, precluding parallelization. Thus, PFIM can leverage parallel computing to achieve further speedup.

\section{Case Study I: Periodic Solutions in Smooth Systems}
\label{section3}
This section uses smooth single-degree-of-freedom (SDOF) systems as illustrative examples to clarify the solution procedure of PFIM, especially for cases with unknown excitation frequencies and continuation along parameter paths. Theoretical convergence properties are examined and validated through numerical results obtained with PCA, highlighting the accuracy and efficiency of the proposed approach in practical applications.

\subsection{Self-Excited Oscillator: Solution and Validation}
\label{section3.1}
A classical van der Pol oscillator is selected as the first benchmark to illustrate the PFIM procedure. The governing equation is written as
\begin{equation}
	\ddot{x} + x + \mu (x^{2}-1)\dot{x} = 0.
	\label{s_1}
\end{equation}
where $\mu > 0$ controls the strength of the nonlinearity and self-excitation, with a value of $0.9$ used here. The objective is to compute the limit cycle and its oscillation frequency using the proposed method.

The PFIM is applied by first assuming an initial periodic function and solving the linearized perturbation equation over one period using the PCA-based matrix exponential integration. The unknown frequency is simultaneously updated according to the phase condition until convergence is achieved. The convergence history of the residual norm is monitored to assess the convergence rate.

As mentioned above, PFIM, employing a first-order perturbation expansion, achieves a theoretical quadratic convergence rate for smooth systems. Specifically, near the exact solution $x_{*}$, the error after one iteration decreases quadratically, i.e., the error at iteration $i$, defined as $e_i = ||x_{*} - x^{(i)}||$, satisfies
\begin{equation}
	\label{new_1}
	e_{i+1} = \mathcal{C} e_i^2 +\mathcal{O}(e_i^2).
\end{equation}

To better quantify this relationship, the solution obtained by HBM is used as the initial guess, while a high-accuracy shooting method solution $X^{SM}$ serves as the reference. The error is then defined as:
\begin{equation}
	 e_{a}(x^{(i)}) = \frac{1}{n_p} \sum_{k=1}^{n_p} \left| x_{k}^{(i)} - x_{k}^{SM} \right|, 
\end{equation}
and for self-excited systems the frequency error is
\begin{equation} 
	e_{\omega}(\omega^{(i)}) = \left| \omega^{(i)} - \omega^{SM} \right|. 
\end{equation}

\begin{figure*}[!htbp]
	\centering
	\begin{subfigure}[t]{0.48\textwidth}
		\centering
		\includegraphics[width=\linewidth]{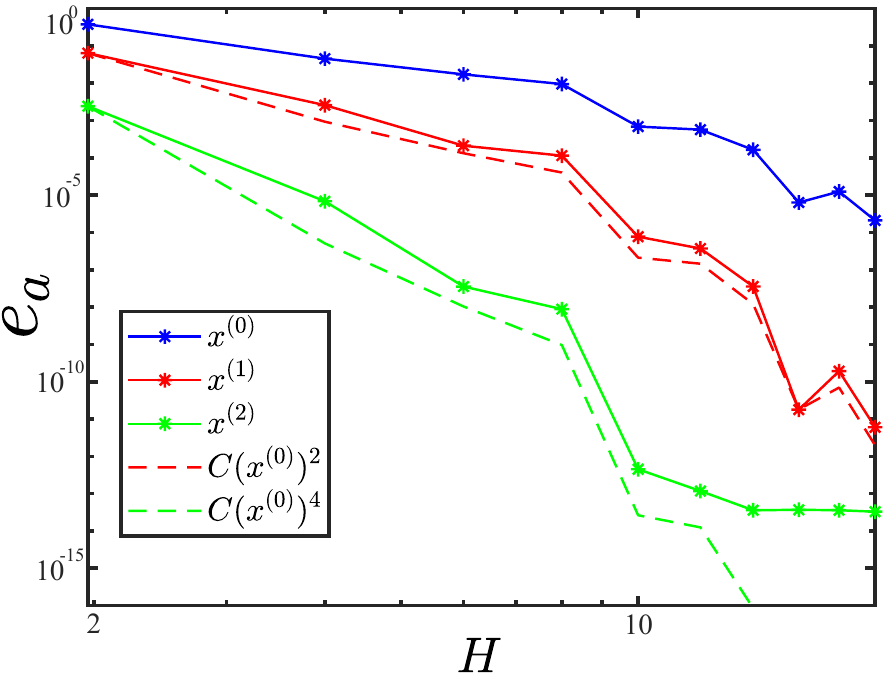}
		\vspace{.2pt}
		\caption{}
		\label{fig.1a}
	\end{subfigure}
	\hfill
	\begin{subfigure}[t]{0.48\textwidth}
		\centering
		\includegraphics[width=\linewidth]{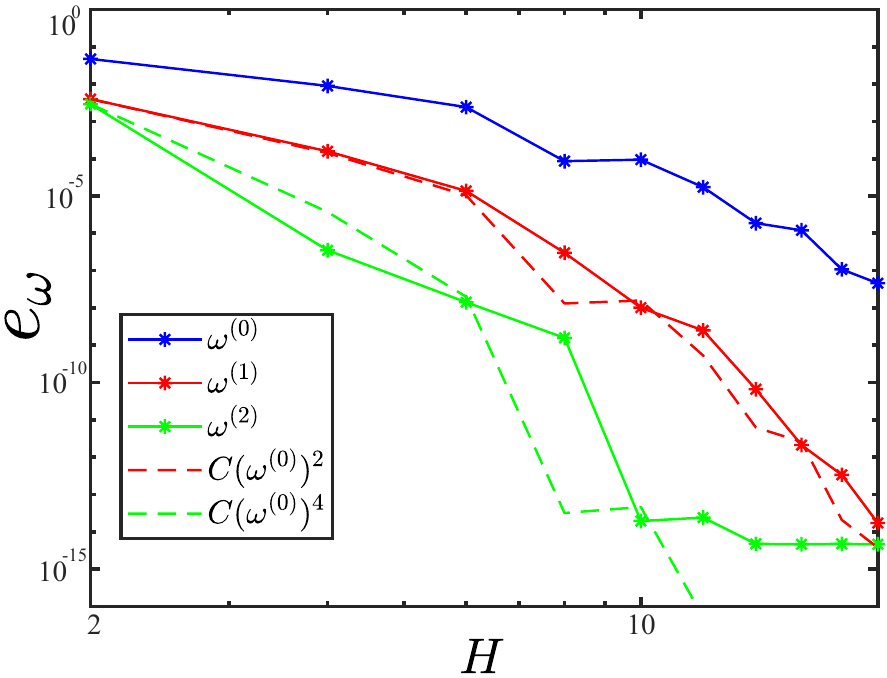}
		\vspace{.2pt}
		\caption{}
		\label{fig.1b}
	\end{subfigure}
	\caption{Convergence of HBM and PFIM solutions for Eq.~\eqref{s_1}. (a) Displacement error: $x^{(0)}$ (HBM); $x^{(1)}$, $x^{(2)}$ (PFIM iterations). Dashed lines $C(x^{(0)})^{2}$, $C(x^{(0)})^{4}$ are scaled error powers. (b) Frequency error with analogous notation.}
	\label{fig_1}
\end{figure*}

The HBM solutions with different orders $H$ are used as the initial guesses $x^{(0)}$ for the PFIM theoretical iteration, and the corresponding iteration errors are plotted in Fig.~\ref{fig_1} (a) and (b). As shown in Fig.~\ref{fig_1} (a), the HBM solutions $x^{(0)}$ exhibit increasingly steeper slopes as $H$ grows, which is consistent with the exponential convergence behavior of HBM for smooth systems. Moreover, excellent agreement is observed between  $x^{(1)}$ and  $C(x^{(0)})^2$, both in their convergence trends and in the magnitude of errors across all harmonic orders. While $x^{(2)}$ remains consistent with  $C(x^{(0)})^4$ for $H \leq 10$, its error plateaus near $10^{-14}$ when $H \leq 10$, suggesting that the iteration has reached the precision limit imposed by the reference solution. Beyond this point, further error reduction becomes unattainable.

Collectively, the results confirm that PFIM attains quadratic convergence for smooth dynamical systems, i.e., 
$e_{i+1} \propto e_{(i)}^{2}$, in agreement with perturbation theory. Having established this theoretical limit, the discussion now turns to the convergence behavior of the numerically implemented PFIM.

In general, for most numerical solution methods such as HBM and SM, the initial guess is typically chosen as the periodic solution of the corresponding linear system. This choice is popular because it is both easy to obtain and reasonably close to the true solution. However, for self-excited vibration systems, the corresponding linear system does not possess a periodic solution. Therefore, in this work, the initial guess is selected as
\begin{equation} 
	x^{(0)} = \cos(\omega^{(0)} t) ,
\end{equation}
where $\omega^{(0)}$ is the corresponding frequency initial guess, and it is fixed at 1 in this example.

\begin{figure*}[!htbp]
	\centering
	\begin{subfigure}[t]{0.48\textwidth}
		\centering
		\includegraphics[width=\linewidth]{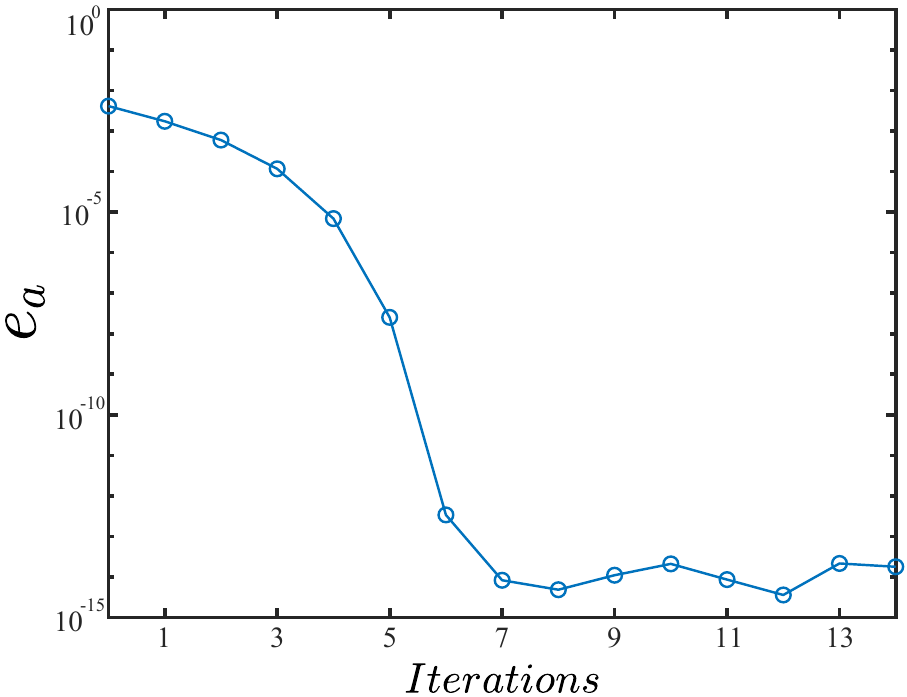}
		\vspace{.2pt}
		\caption{}
		\label{fig.2a}
	\end{subfigure}
	\hfill
	\begin{subfigure}[t]{0.48\textwidth}
		\centering
		\includegraphics[width=\linewidth]{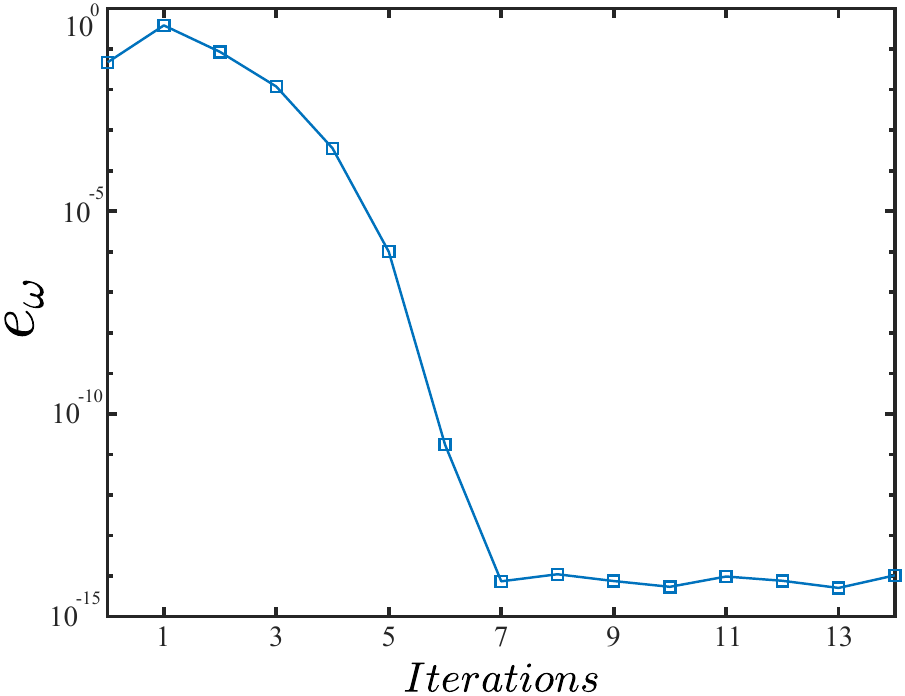}
		\vspace{.2pt}
		\caption{}
		\label{fig.2b}
	\end{subfigure}
	\caption{Error convergence for Eq.~\eqref{s_1} using linear initial guess. (a) Displacement error (b) Frequency error.}
	\label{fig_2}
\end{figure*}

Fig.~\ref{fig_2} illustrates the iterative process when using the PCA-based method under the initial guess described above. It reveals that both displacement and frequency errors attain quadratic convergence rates following the first iteration. By the seventh iteration, these errors reach the precision level of the reference solution, subsequently exhibiting persistent oscillations. These fluctuations arise when the error magnitude approaches machine precision, where further iterations induce computational instabilities due to accumulating floating-point round-off errors.
In this smooth benchmark case, all $n_p$ values are set to $2^{18}$, ensuring that the discussion focuses on the accuracy limit of PFIM. Under these conditions, PFIM demonstrates excellent convergence behavior and accuracy, both in the theoretical results and in the numerical results obtained using the PCA-based method.

\subsection{Parameter Continuation and Frequency Tracking}
\label{section3.2}
Parameter continuation is an important tool for analyzing the dynamic responses of nonlinear systems. A simple SDOF Duffing oscillator is often used to test and benchmark the continuation capability of such methods, since it exhibits several classical nonlinear dynamical features. Its governing equation is given as
\begin{equation}
	\ddot{x} + 0.1\dot{x} + x + 0.1x^{3} = \cos(\omega t).
	\label{s_2}
\end{equation}
where the external excitation frequency $\omega$ is selected as the primary continuation parameter.

\begin{figure*}[!htbp]
	\centering
	\raisebox{0pt}{\begin{subfigure}[t]{0.5\textwidth}
			\centering
			\includegraphics[width=\linewidth]{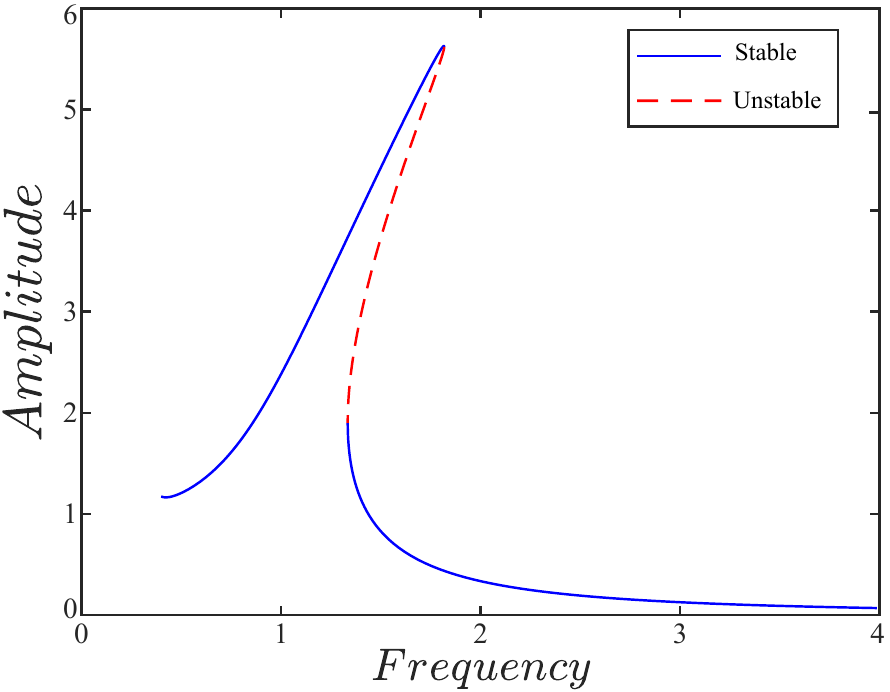}
			\vspace{.2pt}
			\caption{}
			\label{fig.3a}
	\end{subfigure}}
	\hfill
	\begin{subfigure}[t]{0.4125\textwidth} 
		\centering
		\includegraphics[width=\linewidth]{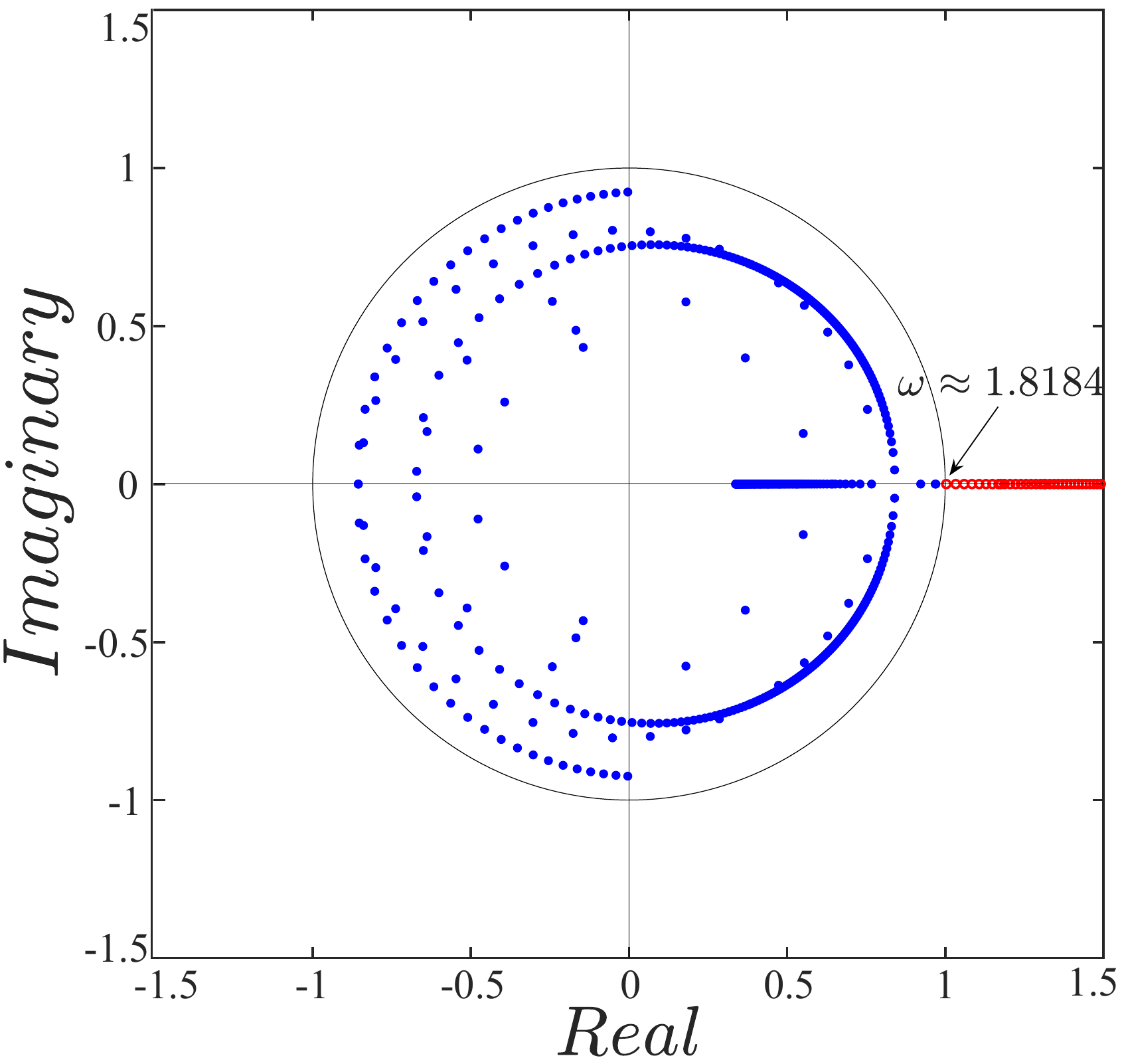}
		\vspace{.2pt}
		\caption{}
		\label{fig.3b}
	\end{subfigure}
	
	\caption{Response of Eq.~\eqref{s_2}. (a) Amplitude-frequency curve (stable: solid, unstable: dashed) (b) Floquet multiplier evolution (filled: stable, hollow: unstable).}
	\label{fig_3}
\end{figure*}

Fig.~\ref{fig_3} (a) and (b) respectively present the amplitude–frequency response curves and the evolution of the corresponding Floquet multipliers for the Duffing system, with the continuation parameter $\omega \in [0.4,4.0]$. In Fig. 3(a), the blue solid lines represent stable periodic solutions, while the red dashed lines indicate unstable ones. The curve exhibits the characteristic hardening stiffness behavior of the Duffing oscillator, including the well-known jump phenomenon in the amplitude response. Transitions between stable and unstable branches correspond to bifurcation events, which are also clearly reflected in Fig.~\ref{fig_3} (b). Specifically, the Floquet multipliers cross $+1$ on the unit circle, indicating the occurrence of a saddle-node (fold) bifurcation in the system.

The total computational cost of PFIM during continuation was further quantified, including overall CPU time, the number of continuation points, and iteration counts. For benchmarking, HBM was selected as a reference method under identical continuation parameters. As previously established, the computational cost of HBM scales approximately as $\mathcal{O}(H^3)$ whereas PFIM scales linearly with $n_p$. To ensure comparable computational complexity, $H^3 \approx n_p$ was imposed as guideline. Accordingly, $n_p = 2^6,2^7,2^8$ were chosen for PFIM, while harmonic numbers $H = 4,5,6$ were selected for HBM.

Tab.~\ref{tab_3} summarizes the continuation results for both methods. From a computational time perspective, neither PFIM nor HBM strictly follows their theoretical complexity scaling at small parameter values. This deviation arises because, when $H$ and $n_p$ are small, the dominant cost is associated with other auxiliary steps rather than the main matrix exponential or linear solves. For PFIM, once $n_p$ increases from $2^7$ to $2^8$, the total computation time nearly doubles, consistent with its linear complexity with respect to the number of sampling points. At low parameter values, the computation times of PFIM and HBM are comparable.

\begin{table}[htbp]
	\centering
	\caption{Computational time, numbers of points on amplitude-frequency response curves, and total numbers of iterations from different methods for the Eq.~\eqref{s_2}}
	\label{tab_3}
	\begin{tabular}{@{}lcccccc@{}}
		\toprule
		& \multicolumn{3}{c}{PFIM} & \multicolumn{3}{c}{HBM} \\
		\cmidrule(lr){2-4} \cmidrule(lr){5-7} 
		 & $n_p=2^6$ & $n_p=2^{7}$ & $n_p=2^{8}$& $H=4$ & $H=5$ & $H=6$ \\
		\midrule
		Computational time  & 0.7148 & 0.9556 & 1.7384 & 0.7571 & 0.8443 & 0.9690\\
		Number of points & 211 & 211 & 211 & 211 & 211 & 211 \\
		Number of iterations & 833 & 721 & 709 & 828 & 825 & 825 \\
		\bottomrule
	\end{tabular}
\end{table}

Regarding the number of continuation points, both methods yield identical results across all parameter settings. This observation confirms that both approaches reliably converge at each continuation step, reflecting their robust convergence performance.

In terms of the total number of iterations, PFIM demonstrates a distinct advantage over HBM, an advantage that becomes more pronounced as the parameters increase. This suggests that PFIM not only maintains comparable accuracy but also achieves greater efficiency in the iterative process, particularly for larger-scale problems.

This section evaluates the continuation capability of PFIM using a SDOF Duffing oscillator as a benchmark example. The results demonstrate that PFIM performs on par with HBM, showing no degradation in accuracy or efficiency. This confirms PFIM’s strong potential for broader applications in smooth nonlinear problems.

\section{Case Study II: Preliminary Application in Non-Smooth Systems}
\label{sec:nonsmooth}

In this section, we further explore the preliminary application of PFIM to non-smooth systems. It is worth emphasizing that PFIM holds significant potential for improvement in this setting — for instance, through techniques such as precise event detection and refined perturbation-based linearization. While previous analyses have rigorously established PFIM’s quadratic convergence in smooth systems — guaranteed by perturbation theory — non-smooth systems exhibit fundamentally different characteristics. The loss of infinite differentiability invalidates Taylor series expansions, rendering direct application of perturbation-based convergence proofs infeasible. This section therefore investigates how system smoothness influences both the theoretical convergence rate and the numerical performance of PFIM, highlighting key challenges and possible pathways for extension.

\subsection{$C^1$-Continuous Systems}
\label{subsec:c1}

Systems with $C^1$-continuity are first considered, where the non-smooth terms in the governing equations are continuously differentiable. As an illustrative example, the term $\dot{x} | \dot{x} |$ models drag forces in viscous fluids: at relatively high velocities, the fluid resistance scales with the square of the velocity and acts in the opposite direction of motion \cite{LAI2024118219}. The governing equation is
\begin{equation}
	\ddot{x} + 0.05\dot{x} + x + 0.5\dot{x} | \dot{x} | = 0.2\cos(t).
	\label{ns_1}
\end{equation}
The partial derivative of the nonlinear term with respect to $\dot{x}$ is 
\begin{equation}
	\frac{\partial}{\partial \dot{x}} \left( \dot{x} | \dot{x} | \right) = 
	\begin{cases} 
		2\dot{x}, & \dot{x} \geq 0 \\
		-2\dot{x}, & \dot{x} < 0 
	\end{cases}
\end{equation}
confirming that the system has continuity $C^1$. PFIM performance is then evaluated using error metrics consistent with those in Figs.~\ref{fig_1} and~\ref{fig_2}. The initial guess for the theoretical solution is taken from the HBM result, whereas the initial guess for the numerical solution is taken from the periodic solution of the underlying linear system.

\begin{figure*}[!htbp]
	\centering
	\raisebox{1pt}{\begin{subfigure}[t]{0.48\textwidth}
			\centering
			\includegraphics[width=\linewidth]{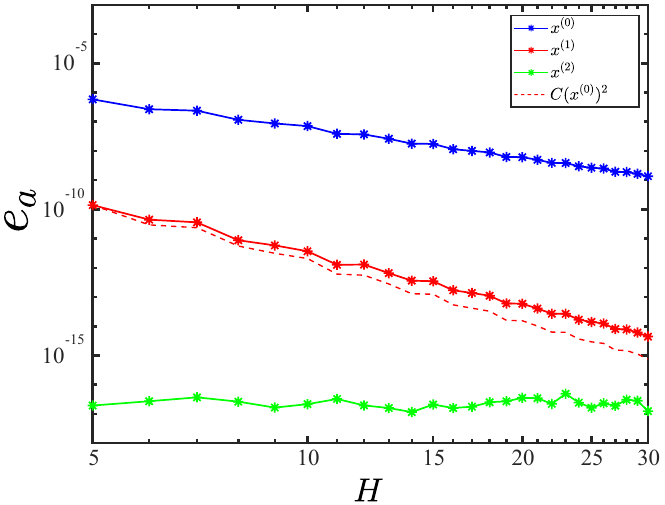}
			\vspace{.2pt}
			\caption{}
			\label{fig.4a}
	\end{subfigure}}
	\hfill
	\begin{subfigure}[t]{0.48\textwidth}
		\centering
		\includegraphics[width=\linewidth]{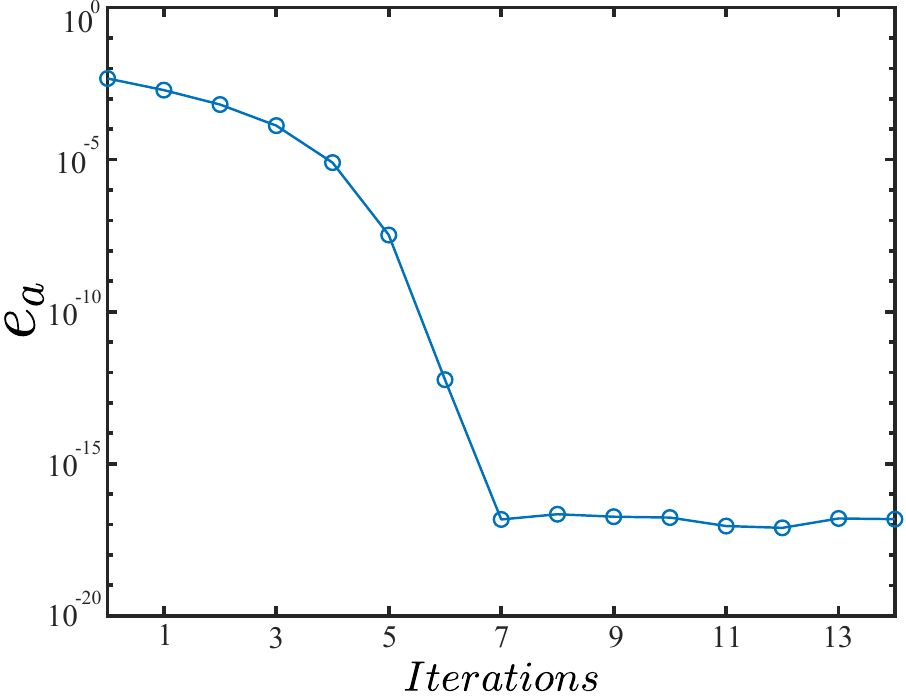}
		\vspace{-0.8pt}
		\caption{}
		\label{fig.4b}
	\end{subfigure}
	
	\caption{Convergence analysis for Eq.~\eqref{ns_1}. (a) HBM ($x^{(0)}$) and PFIM iterations with theoretical solutions ($x^{(1)}$, $x^{(2)}$). Dashed curves indicate scaled error powers; (b) PFIM iterations with Numerical solutions.}
	\label{fig_4}
\end{figure*}

Fig.~\ref{fig_4} (a) shows the evolution of the error during PFIM iterations when HBM solutions of different harmonic orders are used as initial guesses. Detailed convergence behavior of HBM for non-smooth systems can be found in Ref. \cite{wang2021convergence}; here the focus is on the error reduction rate after iteration. As in the smooth case, the curve of $x^{(1)}$ closely follows $C(x^{(0)})^2$, indicating that PFIM retains quadratic convergence for $C^1$-continuous systems. Notably, $x^{(2)}$ reaches machine precision at $H=5$, preventing further comparison with the fourth-order term $C(x^{(0)})^4$.

Fig.~\ref{fig_4} (b) presents numerical results obtained when the linear solution is used as the initial guess. The convergence remains nearly quadratic and stabilizes after approximately seven iterations. The asymptotic error level agrees with theoretical predictions, confirming the consistency between the numerical implementation and the theoretical framework.

Fig.~\ref{fig_4} (a) demonstrates polynomial convergence for HBM ($x^{(0)}$) as harmonic count $H$ increases, consistent with established non-smooth system behavior. The PFIM theoretical solution $x^{(1)}$ closely follows the scaled quadratic term $C(x^{(0)})^{2}$, though with marginally reduced convergence order. This indicates near-but-subquadratic convergence for $C^1$-continuous systems. Notably, $x^{(2)}$ reaches machine precision at $H=5$, precluding comparison with the quartic term $C(x^{(0)})^{4}$.

Taken together, these results demonstrate that for $C^1$-continuous systems, PFIM achieves near-quadratic convergence in both theory and practice. The two approaches reach similar precision limits, with numerical errors approaching machine precision and closely matching the theoretical predictions.

\subsection{$C^0$-Continuous Systems}
\label{subsec:c0}

$C^0$ continuous systems possess lower smoothness than their $C^1$-continuous counterparts,as the non-smooth terms are only continuous but not differentiable, leading to discontinuities in the Jacobian matrix. The term $|x|$ provides a representative example \cite{TIAN2023104478} and is widely used to model piecewise-linear restoring forces. Physically, $|x|$ can describe systems in which the restoring force depends only on the displacement magnitude, irrespective of direction. Typical scenarios include elastic elements with symmetric clearance, unilateral springs, or mechanical stops, where the stiffness abruptly changes when crossing the equilibrium position. Such behavior induces a discontinuous tangent stiffness at $x=0$, which directly affects system dynamics.

\begin{figure*}[!htbp]
	\centering
	\raisebox{1pt}{\begin{subfigure}[t]{0.48\textwidth}
			\centering
			\includegraphics[width=\linewidth]{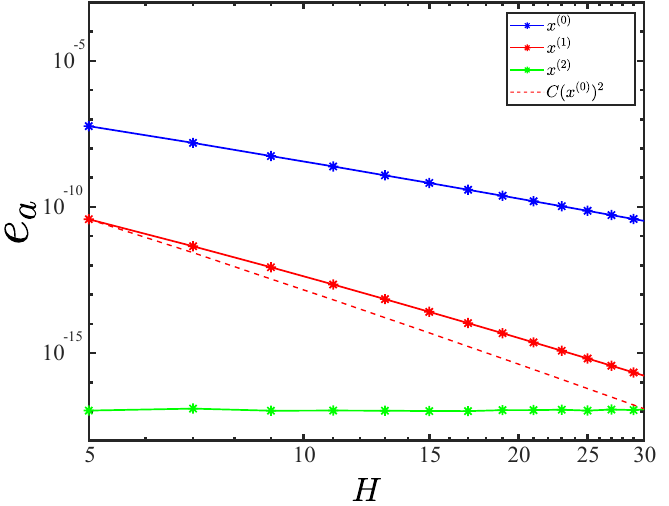}
			\vspace{0pt}
			\caption{}
			\label{fig.5a}
	\end{subfigure}}
	\hfill
	\begin{subfigure}[t]{0.48\textwidth}
		\centering
		\includegraphics[width=\linewidth]{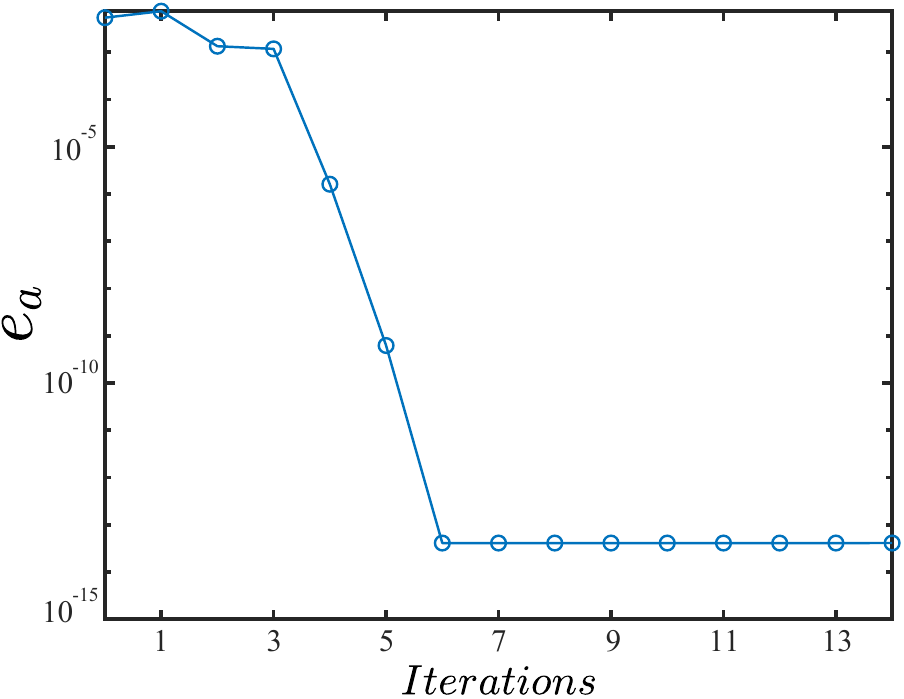}
		\vspace{-0.8pt}
		\caption{}
		\label{fig.5b}
	\end{subfigure}
	
	\caption{Convergence analysis for Eq.~\eqref{ns_2}. (a) HBM ($x^{(0)}$) and PFIM iterations with theoretical solutions ($x^{(1)}$, $x^{(2)}$). Dashed curves indicate scaled error powers; (b) PFIM iterations with numerical solutions.}
	\label{fig_5}
\end{figure*}

Consider the following SDOF system:
\begin{equation}
	\ddot{x} + 0.05\dot{x} + x + 0.5 |x| = 0.2\cos(t).
	\label{ns_2}
\end{equation}
The state-space derivatives of $|x|$ are:
\begin{equation}
	\frac{\partial}{\partial x} |x| = 
	\begin{cases} 
		1, & x > 0 \\
		-1, & x < 0 
	\end{cases}
\end{equation}
exhibiting discontinuity at $x=0$. Fig.~\ref{fig_5} evaluates PFIM performance for this system using the same error metrics as in the previous cases.

Fig.~\ref{fig_5} (a) illustrates the convergence of the theoretical solution for different initial guesses. Unlike the $C^1$-continuous cases, the gap between $x^{(1)}$ and the scaled quadratic curve  $C(x^{(0)})^{2}$ becomes more pronounced. This indicates that the failure of the Jacobian at a finite number of points reduces the PFIM convergence rate from quadratic to superlinear. Nevertheless, this loss of smoothness does not prevent convergence of the numerical solution. As shown in Fig.~\ref{fig_5} (b), PFIM converges within only six iterations when initialized with the periodic solution of the linear system, reaching machine precision. These results confirm that PFIM combined with the PCA integration scheme can effectively handle $C^0$-continuous case. It should be noted, however, that the final attainable precision is slightly lower than in the $C^1$-continuous case, reflecting the influence of neglecting Jacobian discontinuities in the numerical implementation. It is anticipated that developing strategies to compensate for the effect of Jacobian discontinuities at a finite number of points will be an important future research direction for PFIM.

\subsection{$C^{-1}$-Continuous Systems}
\label{subsec:cm1}

Building on the previous subsections, which examined $C^{1}$- and $C^{0}$-continuous systems, we now focus on the most severe case of non-smoothness,namely $C^{-1}$-continuous system. A representative nonlinearity is the sign function $sign(\dot{x})$, which introduces discontinuities directly into the governing equations. In engineering practice, this type of non-smooth nonlinearity commonly arises in systems involving Coulomb friction \cite{DAI2022106932}, dry sliding contacts, relay control laws, and backlash mechanisms. Such systems exhibit discontinuous restoring or dissipative forces that switch instantaneously with the sign of velocity, leading to piecewise-smooth dynamics and potential non-unique trajectories near switching surfaces.

Consider the following single-degree-of-freedom system incorporating this non-smooth nonlinearity:
\begin{equation}
	\ddot{x} + 0.05\dot{x} + x + 0.02\,\text{sign}(\dot{x}) = 0.2\cos(t).
	\label{ns_3}
\end{equation}
The derivative of the nonlinear term $sign(\dot{x})$ vanishes for all $\dot{x} \neq 0$ and is undefined at $\dot{x} = 0$. In practical implementations, the singular point is ignored during Jacobian evaluation. However, this approximation introduces additional uncertainty compared with the $C^{0}$-continuous case, as the system loses information precisely at the velocity-switching surface, where its dynamics are most sensitive.

\begin{figure*}[!htbp]
	\centering
	\raisebox{-0.3pt}{\begin{subfigure}[t]{0.48\textwidth}
			\centering
			\includegraphics[width=\linewidth]{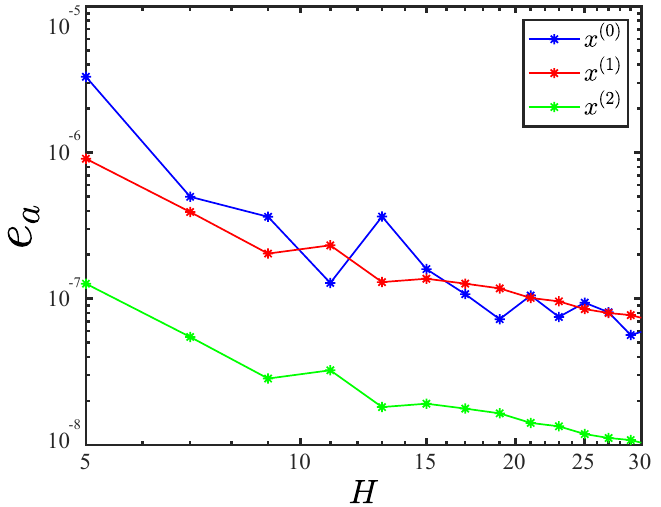}
			\vspace{.2pt}
			\caption{}
			\label{fig.6a}
	\end{subfigure}}
	\hfill
	\begin{subfigure}[t]{0.4785\textwidth}
		\centering
		\includegraphics[width=\linewidth]{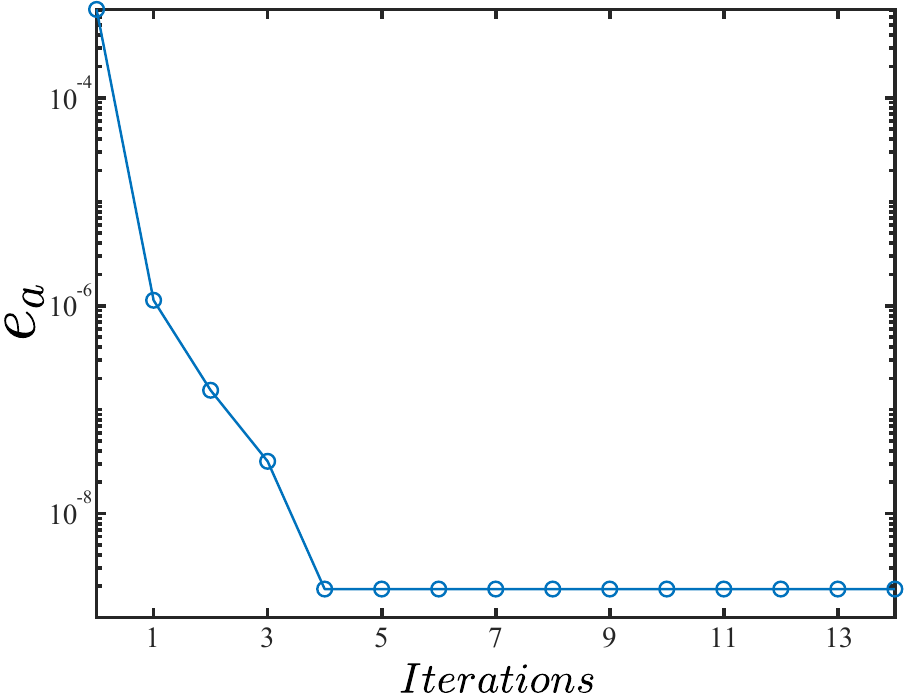}
		\vspace{.5pt}
		\caption{}
		\label{fig.6b}
	\end{subfigure}
	
	\caption{Convergence analysis for Eq.~\eqref{ns_3}. (a) HBM ($x^{(0)}$) and PFIM theoretical solutions ($x^{(1)}$, $x^{(2)}$). Dashed curves indicate scaled error powers; (b) PFIM numerical implementation results.}
	\label{fig_6}
\end{figure*}

Fig.~\ref{fig_6} presents the error convergence of PFIM for both the theoretical and numerical solutions. In Fig.~\ref{fig_6} (a), it is evident that PFIM loses its superlinear convergence and degenerates to linear convergence. In some cases, the error of the first iteration $x^{(1)}$ even exceeds that of the initial guess $x^{(0)}$ obtained from HBM. Moreover, the convergence curves of $x^{(1)}$ and $x^{(2)}$ nearly coincide, further confirming that the convergence rate has degraded to linear.

Fig.~\ref{fig_6} (b) shows the numerical PFIM results, which remain robust and converge within only four iterations. However, consistent with the theoretical prediction, the attainable precision is reduced due to the absence of a well-defined Jacobian at $\dot{x} = 0$. Compared with the $C^{0}$-continuous case, the final error level is further elevated, reflecting the more severe non-smoothness of $C^{1}$-continuous systems.

A more complex $C^{-1}$-continuous systems is further considered, whose governing equations are
\begin{equation}
	\begin{aligned}
		&\ddot{x} + 0.05\dot{x} + x + f_{nl}(x,\dot{x}) = 0.2\cos(3t), \\
		&f_{nl}(x,\dot{x}) = \textstyle\begin{cases}
			2\dot{x}+10x, & G(x,\dot{x}) \geq 0\\
			0, & G(x,\dot{x}) < 0
		\end{cases},
	\end{aligned}
	\label{ns_5}
\end{equation}
where $G(x,\dot{x})=2\dot{x}H(-\dot{x})+10x$, and $H(y)=\{0,y\leq0; 1,y>0 \}$ is the heaviside function \cite{Theodosiou2011}. The nonlinear term will be discontinuous when the function $G(x,\dot{x})=0$. And the partial derivatives of nonlinear term are:
\begin{equation}
	\frac{\partial}{\partial x} f_{nl}(x,\dot{x}) = 
	\begin{cases} 
		10, & G(x,\dot{x}) \geq 0 \\
		0, & G(x,\dot{x}) < 0
	\end{cases},
	\frac{\partial}{\partial \dot{x}} f_{nl}(x,\dot{x}) = 
	\begin{cases} 
		2, & G(x,\dot{x}) \geq 0 \\
		0, & G(x,\dot{x}) < 0
	\end{cases}.
	\label{C-1}
\end{equation}

\begin{figure*}[!htbp]
	\centering
	\raisebox{1pt}{\begin{subfigure}[t]{0.7\textwidth}
			\centering
			\includegraphics[width=\linewidth]{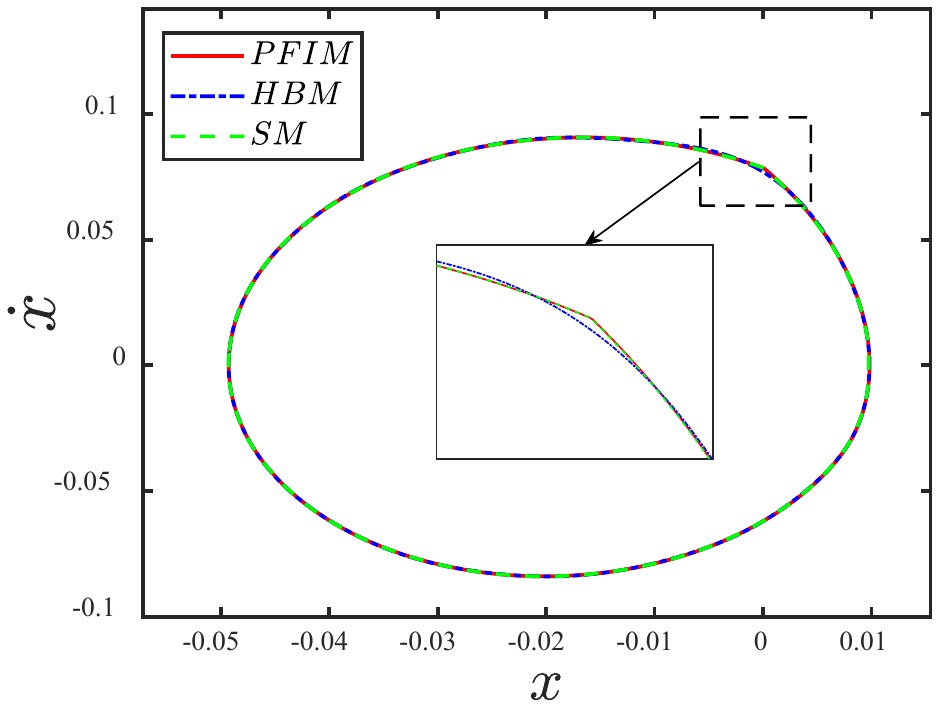}
	\end{subfigure}}	
	\caption{Phase portrait of Eq.~\eqref{C-1}. The red solid line corresponds to PFIM, the blue dash-dotted line corresponds to HBM, and the green dashed line corresponds to SM.}
	\label{fig_C-1_phase}
\end{figure*}

Due to this high degree of discontinuity, the system exhibits pronounced non-smooth behavior, as illustrated in the phase portrait in Fig.~\ref{fig_C-1_phase}. As in the previous examples, a high-precision reference solution is obtained using SM, and HBM employs 10 harmonic terms. The SM solution reveals a distinct discontinuity point in the periodic trajectory (upper-right corner of the phase portrait), which cannot be captured by HBM. In contrast, PFIM shows excellent agreement with the SM results over the entire trajectory, including at the discontinuous point.

As previously noted, HBM faces two inherent and conflicting challenges when applied to non-smooth problems. First, a sufficiently large number of harmonics is required to accurately represent the discontinuous nonlinear restoring force, as illustrated in Fig.~\ref{fig_C-1_phase}. However, increasing the number of harmonics cannot eliminate the Gibbs phenomenon, which is intrinsic to Fourier-series representations of discontinuous functions.

\begin{figure*}[!htbp]
	\centering
	\raisebox{1pt}{\begin{subfigure}[t]{0.65\textwidth}
			\centering
			\includegraphics[width=\linewidth]{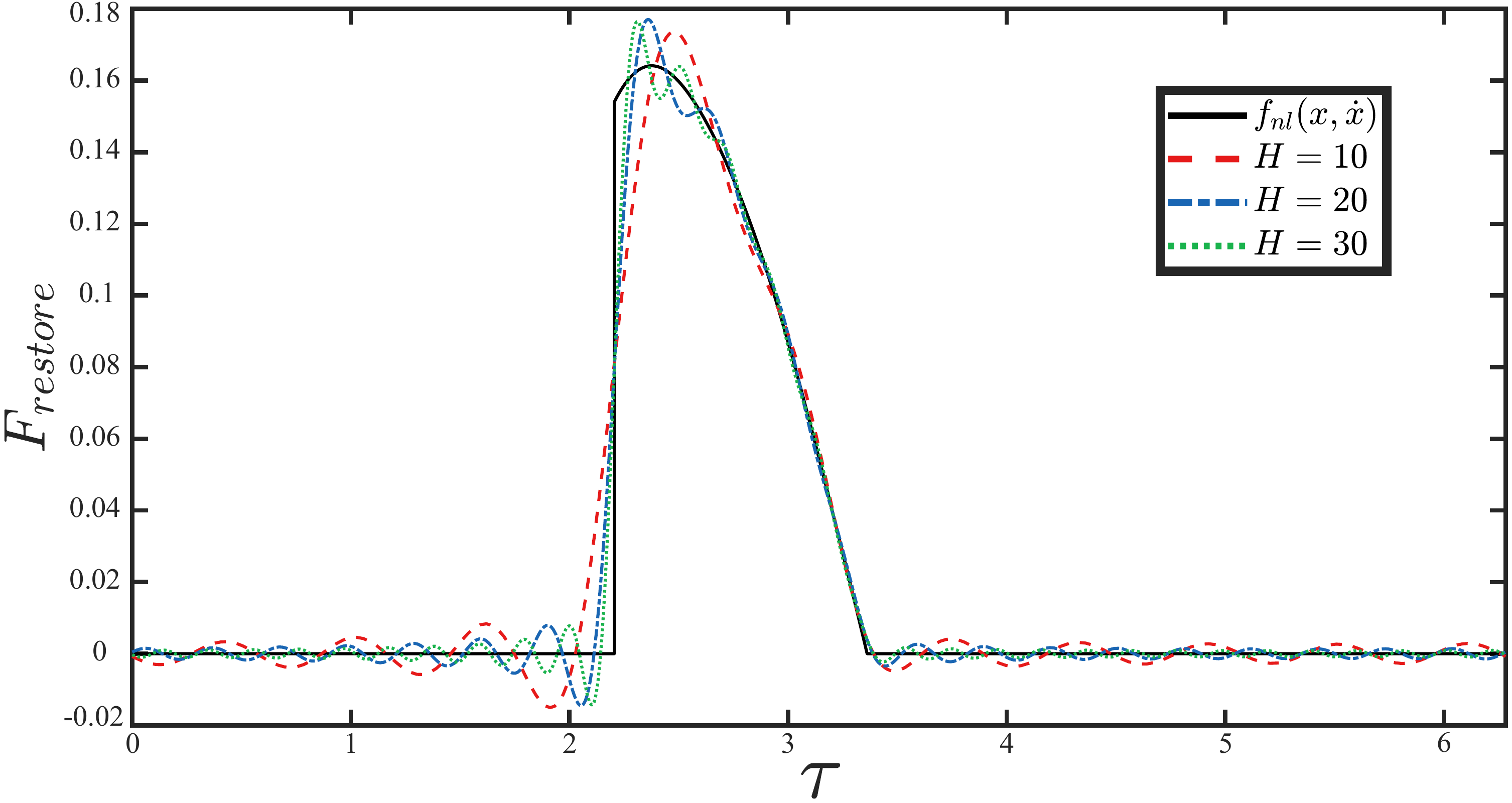}
	\end{subfigure}}	
	\caption{Approximation of the nonlinear restoring force under different harmonic orders. The solid black line represents the target discontinuous nonlinear restoring force. The red dashed curve corresponds to the approximation with $H=10$. The blue dash-dotted curve corresponds to the approximation with $H=20$. The green dotted curve corresponds to the approximation with $H=30$.}
	\label{C-1_nl_force}
\end{figure*}

Figure~\ref{C-1_nl_force} provides a more detailed illustration of this issue by showing the approximation of $f_{nl}$ using $H = 10$, $20$, and $30$ harmonic terms. The black solid line represents the exact non-smooth restoring force, while the red, blue, and green curves correspond to the reconstructed functions obtained from HBM with increasing harmonic orders. Although higher $H$ improves the global approximation, pronounced overshoots and oscillations persist near the discontinuity, characteristic of the Gibbs phenomenon. This overshoot does not vanish with additional harmonics and typically converges to approximately $9\%$ of the jump magnitude, which can lead to significant local errors in the predicted response \cite{gottlieb2011review}.

\begin{figure*}[!htbp]
	\centering
	\raisebox{2pt}{\begin{subfigure}[t]{0.48\textwidth}
			\centering
			\includegraphics[width=\linewidth]{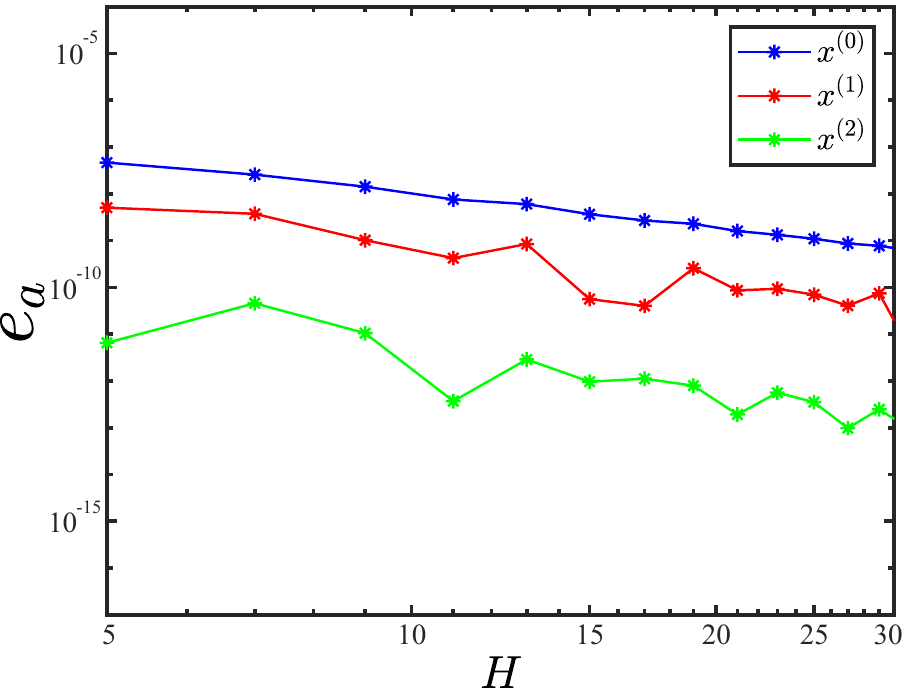}
			\vspace{.25pt}
			\caption{}
			\label{C-1_x_a}
	\end{subfigure}}
	\hfill
	\begin{subfigure}[t]{0.46\textwidth}
		\centering
		\includegraphics[width=\linewidth]{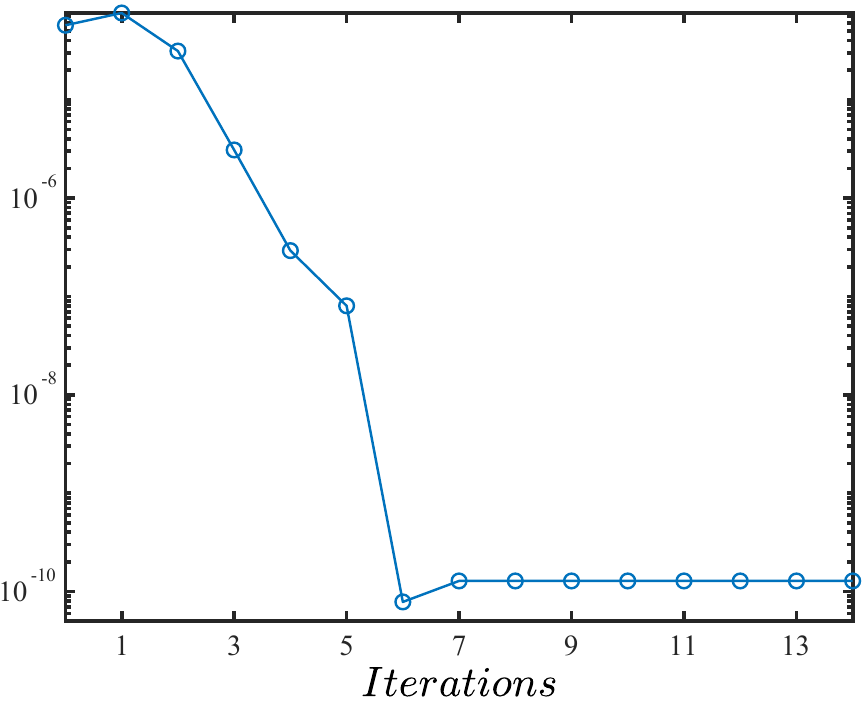}
		\vspace{-2pt}
		\caption{}
		\label{C-1_x_b}
	\end{subfigure}
	
	\caption{Convergence analysis for Eq.~\eqref{C-1}. (a) HBM ($x^{(0)}$) and PFIM theoretical solutions ($x^{(1)}$, $x^{(2)}$). Dashed curves indicate scaled error powers; (b) PFIM numerical implementation results.}
	\label{C-1_x}
\end{figure*}

Compared with HBM, PFIM captures the non-smooth discontinuities in the response through its time-domain formulation and basis-free numerical scheme. Nevertheless, the solution process still relies on the Jacobian matrix and the continuity of the response. Figure~\ref{C-1_x} illustrates the convergence of PFIM for both the theoretical and numerical solutions. Similar to the previous $C^{-1}$-continuous system, PFIM exhibits a linear convergence rate. Moreover, the attainable precision of the numerical solution is further reduced due to the more severe derivative discontinuities.

In addition to the non-smooth functions of $x$ and $\dot{x}$ considered above, another important class of non-smooth problems involves time-dependent discontinuities. A typical example is the externally forced system subjected to a square-wave excitation, 
\begin{equation}
	\ddot{x} + 0.05\dot{x} + x + 0.1x^3 = 0.2F(t), \quad F(t) = 
	\begin{cases} 
		1, & 0 \leq t < T/2 \\
		0, & T/2 \leq t < T 
	\end{cases}
	\label{ns_4}
\end{equation}
where $T=4\pi$. Square-wave forcing is widely used to model systems subjected to intermittent or on–off excitations, such as relay-controlled actuators, pulse-width-modulated loads, and impact-driven mechanisms. It is worth noting that this type of non-smooth term has no effect on the Jacobian matrix, as it does not involve the state variables $x$ or $\dot{x}$.

\begin{figure*}[!htbp]
	\centering
	\raisebox{0pt}{\begin{subfigure}[t]{0.475\textwidth}
			\centering
			\includegraphics[width=\linewidth]{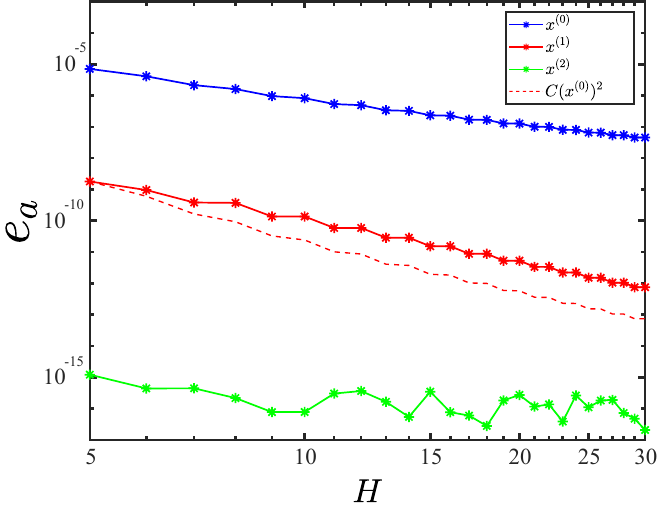}
			\vspace{-5mm}
			\caption{}
			\label{fig.7a}
	\end{subfigure}}
	\hfill
	\begin{subfigure}[t]{0.475\textwidth}
		\centering
		\includegraphics[width=\linewidth]{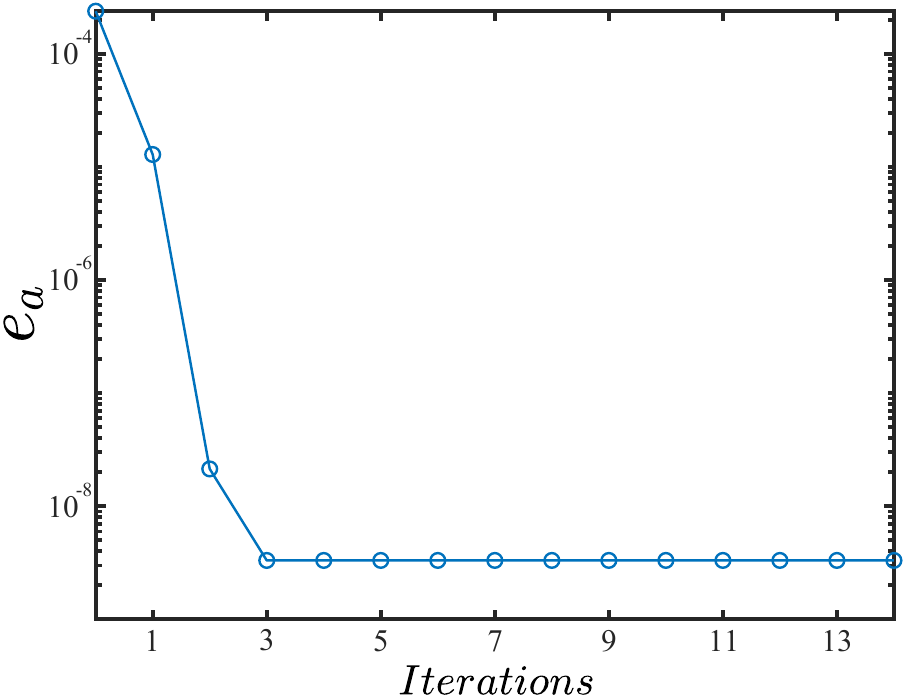}
		\caption{}
		\label{fig.7b}
	\end{subfigure}
	
	\caption{Convergence analysis for Eq.~\eqref{ns_4}. (a) HBM ($x^{(0)}$) and PFIM theoretical solutions ($x^{(1)}$, $x^{(2)}$).  Dashed curves indicate scaled error powers; (b) PFIM numerical implementation results.}
	\label{fig_7}
\end{figure*}
As shown in Fig.~\ref{fig_7} (a), the non-smooth external excitation does not degrade the convergence rate of the theoretical PFIM solution as severely as the $C^{-1}$-continuous functions of the state variables. PFIM retains superlinear convergence, approaching second order. This behavior is closely related to the fact that the square-wave excitation does not affect the Jacobian matrix.

However, the strong non-smoothness of the excitation still limits the attainable precision of the numerical solution. As illustrated in Fig.~\ref{fig_7} (b), the final error level is comparable to that observed in other $C^{-1}$-continuous systems, indicating that the numerical implementation remains sensitive to discontinuities in the excitation.

In summary, the preliminary application of PFIM to non-smooth systems allows several fundamental numerical observations to be drawn. The convergence behavior of PFIM is governed by the smoothness of the nonlinear terms. Specifically, PFIM exhibits quadratic convergence for smooth systems and retains this rate for $C^{1}$-continuous nonlinearities involving state variables. As the degree of non-smoothness increases, such as in $C^{0}$- and $C^{-1}$-continuous systems, the convergence rate degrades to superlinear or even linear.

In addition to the convergence rate, the attainable precision of the solution is also affected by the presence of non-smooth terms. Regardless of whether the non-smoothness is associated with state variables or external excitation, it consistently limits the final accuracy of the solution, which distinguishes its effect from that on the convergence rate.

\section{Case Study III: High-Dimensional Finite-Element Example}
\label{sec:highdim}

\begin{figure*}[!htbp]
	\centering
	\includegraphics[width=0.65\textwidth]{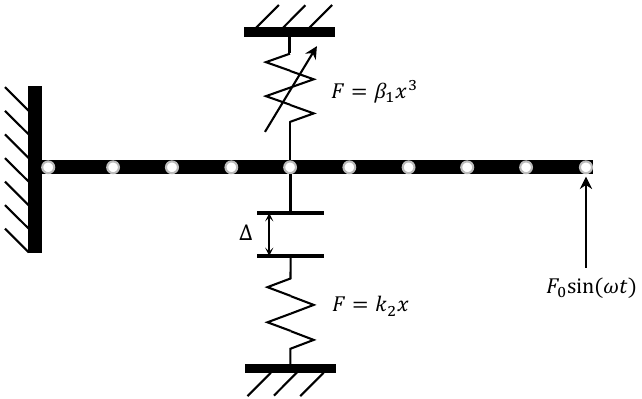}
	\caption{Finite element model of a cantilever beam with 10 nodes. A cubic spring and gap spring act at node 5, while sinusoidal excitation is applied at node 10.}
	\label{fig_8}
\end{figure*}

Finally, a high-dimensional finite element example is employed to demonstrate the potential of PFIM in practical engineering applications. As illustrated in Fig. \ref{fig_8}, the structure is fixed at the left end and consists of ten nodes. Harmonic excitation is applied in the vertical direction at node 10. The nonlinearities arise from a cubic spring and a clearance spring attached to node 5. The cubic nonlinearity can represent the geometric nonlinearity of an inclined cable or similar structural components, while the clearance spring models effects such as assembly imperfections or base motion. The governing equation is as follows:
\begin{equation}
	\mathbf{M\ddot{X}} + \mathbf{C\dot{X}} +\mathbf{KX} + \mathbf{F}_{ns}(\mathbf{X}) = \mathbf{F}_e
	\label{hs}
\end{equation}
where $\mathbf{X} = [y_2, \theta_2, \dots, y_{10}, \theta_{10}]^{\mathrm{T}}$ contains vertical displacements $y_i$ (positive upward) and rotation angles $\theta_i$ (positive counterclockwise) for nodes 2-10 (node 1 is fixed). The nonlinear force $\mathbf{F}_{ns}$ acts only at node 5:
\begin{equation}
	F_{ns5} = 
	\begin{cases}
		\beta_{1}y_5^3 + k_2(y_5 + \delta), & y_5 < -\delta \\
		\beta_{1}y_5^3, & y_5 \geq -\delta
	\end{cases}
\end{equation}
with all other zero elements. An external harmonic excitation is applied at the free end of the structure, 
\begin{equation}
	F_{e17} = F_0 \sin(\omega t),
\end{equation}
while all other degrees of freedom remain unforced. The key parameters are: $F_0 = 100$ N, $\omega = 1$ rad/s, $\beta_1 = 1 \times 10^6$ N/m$^3$, $k_2 = 5 \times 10^3$ N/m, $\delta = 0.01$ m. The mass matrix $\mathbf{M}$, the damping matrix $\mathbf{C}$, and the stiffness matrix $\mathbf{K}$ are defined in the Appendix~\ref{Appendix}.

\begin{figure*}[!htbp]
	\centering
	\begin{subfigure}[t]{0.475\textwidth}
			\centering
			\includegraphics[width=\linewidth]{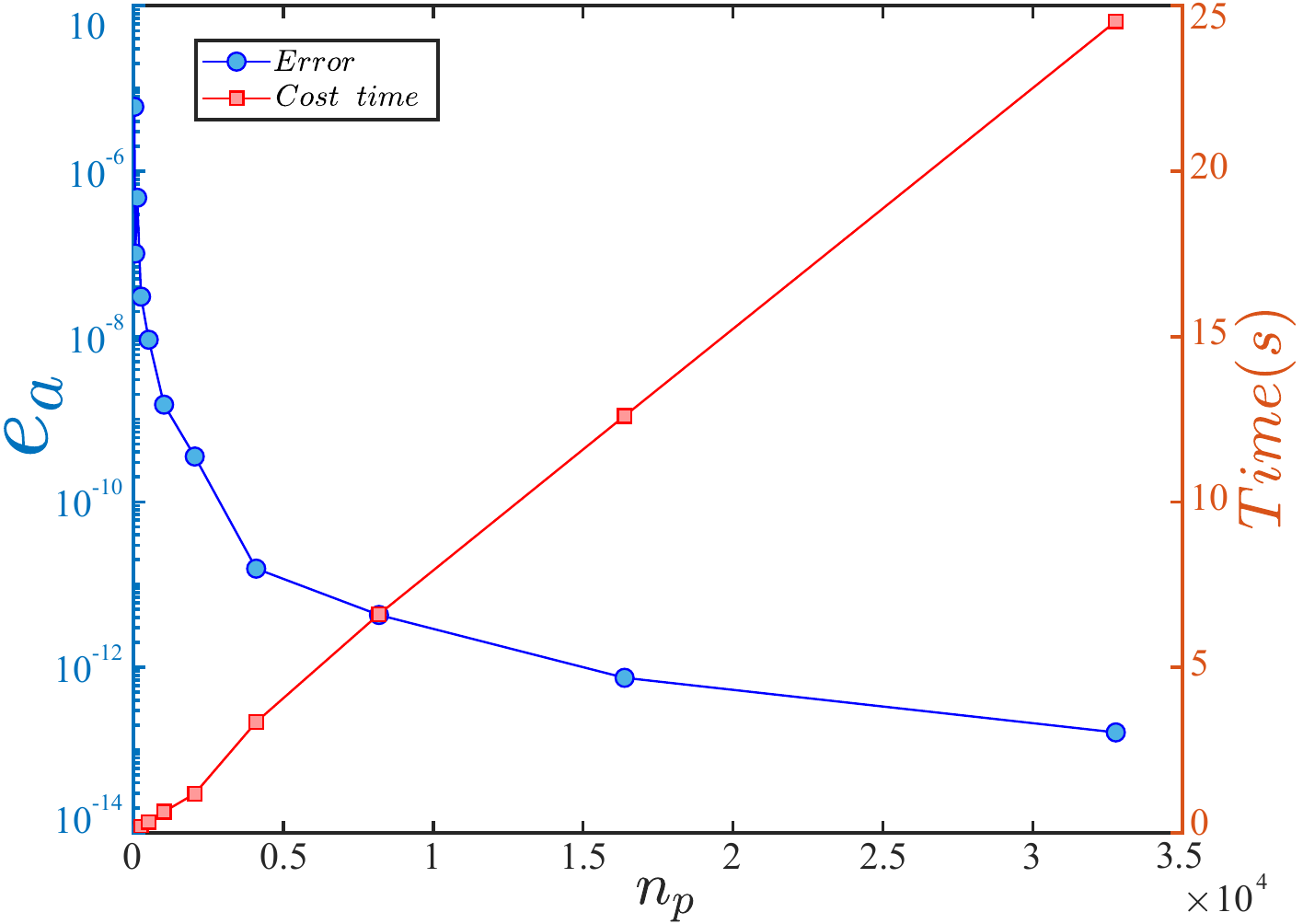}
			\vspace{0pt}
			\caption{}
			\label{fig.9a}
	\end{subfigure}
	\hfill
	\raisebox{-2pt}{\begin{subfigure}[t]{0.48\textwidth}
		\centering
		\includegraphics[width=\linewidth]{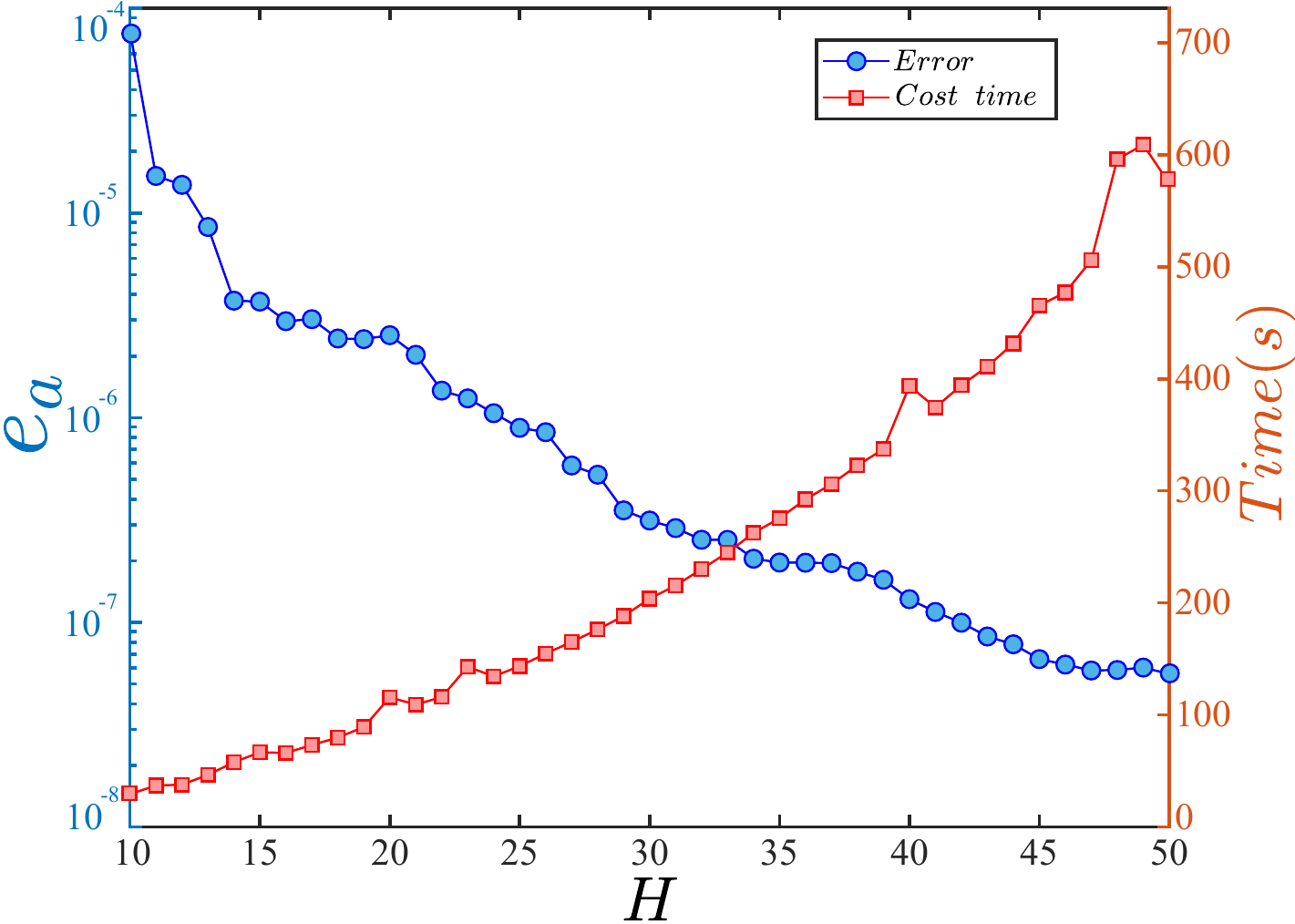}
		\vspace{.2pt}
		\caption{}
		\label{fig.9b}
	\end{subfigure}}
	
	\caption{Convergence analysis for Eq.~\eqref{hs}. (a) PFIM: error and time versus iterations; (b)HBM: error (left axis) and computation time (right axis) versus harmonics.}
	\label{fig_9}
\end{figure*}

To better demonstrate the efficiency advantage of PFIM over HBM, Fig.~\ref{fig_9} presents the convergent results of both PFIM and HBM as $n_p$ and $H$ increase, respectively. The left vertical axes in subfigures (a) and (b) represent the accuracy difference between the final converged result and the result obtained by the high-precision SM method, while the right vertical axes show the computational time required.

From Fig.~\ref{fig_9} (a), it can be clearly observed that as $n_p$ increases from $2^5$ to $2^{15}$, the final accuracy of PFIM gradually stabilizes (approximately $\mathcal{O}(10^{-12})$), and its required time increases proportionally with $n_p$, up to approximately 25 seconds. In contrast, Fig.~\ref{fig_9} (b) shows that for HBM, as $H$ increases from 10 to 50, the computational time increases cubically (peaking at about 600 seconds), and the achieved accuracy only reaches a minimum of $\mathcal{O}(10^{-7})$.

It is worth mentioning that the curves depicting time increase versus $n_p$ or $H$ for PFIM and HBM in Fig.~\ref{fig_9} do not perfectly align with the computational complexity relationships discussed in Section 2. This discrepancy arises because the number of iterations required for convergence may vary for PFIM and HBM under different $n_p$ or $H$. However, the results of both methods generally correspond to our conclusions regarding computational complexity. In summary, the results in Figure 12 strongly confirm the computational advantage of PFIM over HBM in high-dimensional non-smooth systems. Comparing the two subfigures reveals that PFIM requires approximately two orders of magnitude less time than HBM to achieve a similar order of accuracy.

An additional note is that the HBM in this example uses a numerical Jacobian matrix instead of an analytical one, which imposes a significant computational burden. However, obtaining the analytical Jacobian for HBM in non-smooth systems is typically challenging, as it requires precisely locating the piecewise points and performing piecewise integration across different intervals during the Galerkin process. Therefore, for generality, HBM employed the numerical Jacobian in this computation. In contrast, the analytical Jacobian for PFIM can be easily derived simply by taking partial derivatives of the original equations. Thus, it can be argued that PFIM holds a clear computational efficiency advantage over HBM partly because its analytical Jacobian is more readily available.

\begin{figure*}[!htbp]
	\centering
	\includegraphics[width=0.65\textwidth]{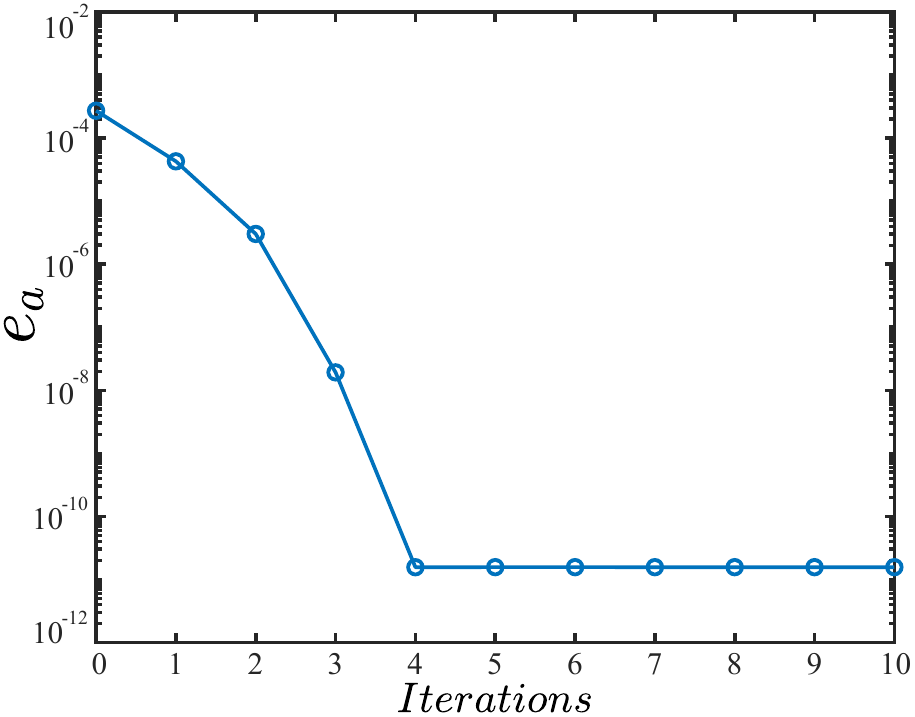}
	\caption{Convergence analysis of the PFIM numerical results for Eq.~\eqref{hs} with $n_p = 2^{12}$.}
	\label{fig_10}
\end{figure*}

Furthermore, we examine the convergence capability and continuation capability of PFIM in this system. Since HBM exhibited very slow convergence when using the solution of the underlying linear system as initial guess, the results in Fig.~\ref{fig_9} were actually obtained using the HBM solution at $H=5$ as initial guess. In contrast, Fig.~\ref{fig_10} demonstrates the convergence of PFIM with $n_p=2^{12}$ under the initial condition of the underlying linear system solution. PFIM converged in merely 4 iterations, with the error decreasing monotonically throughout the process. This clearly illustrates the superior convergence performance of PFIM compared to HBM in this example.

\begin{figure*}[!htbp]
	\centering
	\includegraphics[width=0.65\textwidth]{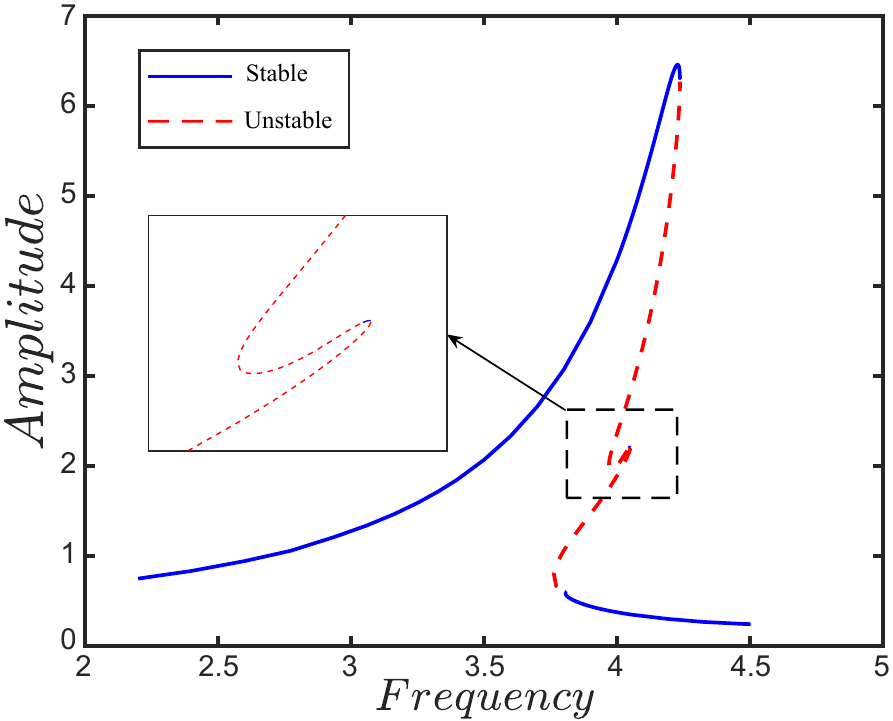}
	\caption{Amplitude-frequency curve of Eq.~\eqref{hs} (stable: solid, unstable: dashed).}
	\label{fig_11}
\end{figure*}

Subsequently, Fig.~\ref{fig_11} presents the continuation curve tracked by PFIM using the external excitation frequency $\omega$ (increasing from 2.2 to 4.5) as the continuation parameter. Consistent with the previous Duffing example, the blue solid line represents stable periodic responses, and the red dashed line represents unstable ones. The continuation curve here is more complex than that of the Duffing system: after reaching the peak (around $\omega \approx 4.2$), the periodic response loses stability and begins to backtrack. When $\omega$ decreases to approximately 3.98, the response amplitude then increases again with further increase in $\omega$. Furthermore, the zoomed-in view shows that during the ascent to a local peak, the periodic response remains stable within a very narrow interval before losing stability again and continuing to backtrack, behavior consistent with the subsequent pattern observed in the Duffing system. The successful tracking of this complex continuation curve demonstrates PFIM's robust continuation capability even in high-dimensional non-smooth systems with complex responses.

In conclusion, based on the above results, PFIM demonstrates excellent performance in the high-dimensional non-smooth finite element example. Compared to HBM, it not only offers higher computational efficiency and better convergence but also successfully tracks the complex continuation curve of the system. These results fully illustrate the solving advantages of PFIM in practical engineering applications.

\section{Conclusions}

This study introduces a novel method, termed the Perturbation Function Iteration Method (PFIM), designed for the efficient and accurate computation of periodic responses in nonlinear and even non-smooth dynamical systems. Compared with the Harmonic Balance Method (HBM) and the Shooting Method (SM), PFIM offers three distinct methodological advantages.
First and foremost, the Jacobian matrix in PFIM can be directly obtained by differentiating the original governing equations, thereby eliminating the need to compute the Jacobian with respect to auxiliary algebraic variables such as harmonic coefficients or initial conditions—a process that is often complex and computationally demanding.
Second, PFIM requires basis-free function assumptions in constructing periodic solutions. This basis-free formulation fundamentally avoids issues such as the Gibbs phenomenon, which typically hinders the convergence of Fourier-based methods like HBM when applied to non-smooth systems.
Finally, the core iteration procedure of PFIM involves only the solution of linear systems, which, when combined with the Piecewise Constant Approximation (PCA), yields substantially higher computational efficiency compared with SM.

Based on the proposed framework, a series of representative numerical examples were conducted to verify the efficiency and accuracy of PFIM, leading to several key numerical findings:
\begin{itemize}

\item In smooth systems, PFIM demonstrates quadratic convergence, with numerical results matching the theoretical accuracy limits, and exhibits excellent continuation capability comparable to HBM.

\item The convergence behavior of PFIM is strongly influenced by the smoothness of the system. For $C^{1}$ systems, it maintains quadratic convergence. For $C^{0}$ systems, the method shows superlinear but sub-quadratic convergence. For general $C^{-1}$ systems, the convergence degrades to linear. Crucially, in all tested non-smooth cases, PFIM accurately captures the discontinuities in the system response—an aspect where HBM suffers from significant Gibbs oscillations.

\item PFIM demonstrates outstanding computational performance in large-scale, non-smooth finite element problems. To achieve the same solution accuracy, HBM requires up to two orders of magnitude more computational time than PFIM, primarily due to its reliance on numerically approximated Jacobian matrices. Furthermore, PFIM successfully traces bifurcation branches in such complex systems, underscoring its practical applicability to high-dimensional non-smooth dynamics.

\end{itemize}

Beyond its demonstrated numerical performance, PFIM possesses strong potential for further development. Its algorithmic structure supports model-order reduction and lends itself naturally to parallelization due to its piecewise-based formulation.

In summary, PFIM provides excellent performance in both smooth and non-smooth dynamical systems. Its combination of high efficiency, superior convergence, and robust continuation capability establishes it as a powerful and competitive tool for tackling the challenges of high-dimensional and non-smooth dynamics in modern engineering applications. Future research will focus on enhancing its accuracy in non-smooth problems and realizing efficient large-scale parallel computation.

\section*{Appendix: Finite Element Matrices}
\label{Appendix}
The global matrices $\mathbf{M}$, $\mathbf{C}$, and $\mathbf{K}$ for the cantilever beam model are defined as follows:

\subsection*{Beam Parameters}
\begin{itemize}
	\item Length $L = 8$ m, cross-section $b \times h = 0.02 \times 0.2$ m
	\item Material properties: $E = 3 \times 10^9$ Pa, $\rho = 7800$ kg/m$^3$
	\item Damping ratio $\zeta = 0.02$
	\item Element length $L_e = L/9 = 8/9$ m
\end{itemize}

\subsection*{Element Matrices}
Element stiffness matrix $\mathbf{k}_e$ and mass matrix $\mathbf{m}_e$:
\begin{align*}
	\mathbf{k}_e &= \frac{EI}{L_e^3}
	\begin{bmatrix}
		12 & 6L_e & -12 & 6L_e \\
		6L_e & 4L_e^2 & -6L_e & 2L_e^2 \\
		-12 & -6L_e & 12 & -6L_e \\
		6L_e & 2L_e^2 & -6L_e & 4L_e^2
	\end{bmatrix} \\
	\mathbf{m}_e &= \frac{\rho A L_e}{420}
	\begin{bmatrix}
		156 & 22L_e & 54 & -13L_e \\
		22L_e & 4L_e^2 & 13L_e & -3L_e^2 \\
		54 & 13L_e & 156 & -22L_e \\
		-13L_e & -3L_e^2 & -22L_e & 4L_e^2
	\end{bmatrix}
\end{align*}
where $A = b \times h = 4 \times 10^{-3}$ m$^2$, $I = \frac{b h^3}{12} = 1.333 \times 10^{-5}$ m$^4$.

\subsection*{Global Matrix Assembly}
The global mass matrix $\mathbf{M}$ and stiffness matrix $\mathbf{K}$ ($18 \times 18$) are assembled from 9 identical elements:
\begin{align*}
	\mathbf{M} &= \bigoplus_{e=1}^{9} \mathbf{m}_e \\
	\mathbf{K} &= \bigoplus_{e=1}^{9} \mathbf{k}_e
\end{align*}

\subsection*{Damping Matrix}
The damping matrix uses Rayleigh damping:
\begin{equation*}
	\mathbf{C} = \alpha \mathbf{M} + \beta \mathbf{K}
\end{equation*}
where $\alpha = 0.362$, $\beta = 5.23 \times 10^{-4}$ are calculated from:
\begin{align*}
	\alpha &= \frac{2\zeta\omega_1\omega_2}{\omega_1 + \omega_2}, \quad 
	\beta = \frac{2\zeta}{\omega_1 + \omega_2} \\
	\omega_1 &= 10.49\ \text{rad/s}, \quad \omega_2 = 65.97\ \text{rad/s}
\end{align*}

\subsection*{Degree of Freedom Mapping}
\begin{align*}
	& \text{Node } i: \text{DOF } (2i-1) = y_i,\ (2i) = \theta_i \\
	& \text{Node 5: DOF 7} = y_5 \\
	& \text{Node 10: DOF 17} = y_{10}
\end{align*}

\section*{Acknowledgement}
Supported by National Natural Science Foundation of China (No. 11572356, 11702336, No. 12172387, No. 12572033), Guangdong Province Natural Science Foundation (No. 2016A020223006), and Fundamental Research Funds for Central Universities (No. 17lgjc42, 17lgpy54).

\section*{Data Availability Statements}
Data sharing not applicable to this article as no datasets were generated or analysed during the current study.

\section*{Conflicts of interest}
The authors declare that they have no conflict of interest.

\end{document}